\documentclass[11pt]{article}
 
\usepackage{graphicx}
\usepackage{amssymb, amsmath,}
                                                                                                                                                                                                                                                                                                                                                                                                                                                                                                                                                                                                                                                                                                                                                                                                                                                                                                                                                                                                                                                                                                                                               \newtheorem{theorem}{Theorem}[section]
\newtheorem{corollary}[theorem]{Corollary}
\newtheorem{definition}[theorem]{Definition}
\newtheorem{lemma}[theorem]{Lemma}
\newtheorem{Remark}[theorem]{Remark}
\newtheorem{Proposition}[theorem]{Proposition}
\title{\textbf{Codes defined by forms of degree 2 on quadric and non-degenerate hermitian varieties in $\mathbb{P}^{4}(\mathbb{F}_q)$}}          
\date{}
\author{\textbf{Fr\'ed\'eric A. B. Edoukou}  \\ \\ CNRS, Institut de Math\'ematiques de Luminy  \\ Luminy case 907  - 13288  Marseille Cedex 9 - France \\  E.mail : edoukou@iml.univ-mrs.fr }
\begin{document}
\maketitle

{\footnotesize \begin{flushleft}
\textbf{Abstract}
\end{flushleft}
We study the functional codes of second order defined by G. Lachaud on $\mathcal{X} \subset {\mathbb{P}}^4(\mathbb{F}_q)$ a quadric of rank($\mathcal{X}$)=3,4,5 or a non-degenerate hermitian variety. We give some bounds for 
the number of points of quadratic sections of $\mathcal{X}$, which are the best possible  and show that codes defined on non-degenerate quadrics are better than those defined on degenerate quadrics. We also show the geometric structure of the minimum weight codewords and estimate the second weight of these codes. For $\mathcal{X}$ a non-degenerate hermitian variety, we list the first five weights and the corresponding codewords. The paper ends with two conjectures. One on the minimum distance for the functional codes of order $h$ on $\mathcal{X} \subset {\mathbb{P}}^4(\mathbb{F}_q)$ a non-singular hermitian variety. The second conjecture on the distribution of the codewords of the first five weights of the  functional codes of second order on $\mathcal{X} \subset {\mathbb{P}}^N(\mathbb{F}_q)$ the non-singular hermitian variety.
}
\section{Introduction}
We study the codes $C_2(\mathcal{X})$ defined on the quadric $\mathcal{X}$ and on the non-singular hermitian variety  $\mathcal{X}: x_0^{t+1}+x_1^{t+1}+x_2^{t+1}+x_3^{t+1}+x_4^{t+1}=0$ in 
$\mathbb{P}^4(\mathbb{F}_q)$ with $q=t^2$ ($t$ is a prime power). \\
When $\mathcal{X}$ is a quadric, we will consider the cases where rank($\mathcal{X}$)= 3, 4, 5, since otherwise, the computation of the miminum distance of the codes $C_2(\mathcal{X})$ is reduced to the study of plane or hyperplane sections of quadrics which is well known.\\
The case where $\mathcal{X}$ is a non-degenerate quadric in $\mathbb{P}^4(\mathbb{F}_q)$ (i.e. rank$(\mathcal{X})=5$, $\mathcal{X}$ is a parabolic quadric), a bound for the minimum distance is given by D. B. Leep and L. M. Schueller  \lbrack 11, p.172\rbrack : 
$$\# \mathcal{X}_{ Z(\mathcal{Q})}(\mathbb{F}_{q} ) \le 3q^2+q+1.$$ 
It is not optimal. When $\mathcal{X}$ is a degenerate quadric of rank 3 or 4, the hypothesis of D. B. Leep and L. M. Schueller is too much restrictive and we cannot find the minimal distance. In this case, however, we have a bound given by G. Lachaud \lbrack 10, proposition 2.3 \rbrack: $$\# \mathcal{X}_{ Z(\mathcal{Q})}(\mathbb{F}_{q} ) \le 4(q^2+q+1)$$ which is not the best possible.\\ 
In the case where $\mathcal{X}$ is a non-singular hermitian variety, the result of Lachaud gives an upper bound \lbrack 14, p.208\rbrack, $$\# \mathcal{X}_{ Z(\mathcal{Q})}(\mathbb{F}_{q} ) \le 2(t+1)(q^2+q+1)$$ which is also not optimal.\\
In this paper we find some optimal bounds excepting the case where $\mathcal{X}$ is a degenerate quadric with rank$(\mathcal{X})=4$ and projective index g$(\mathcal{X})$=1.\\
The paper is organized as follows. First of all we recall some generalities on quadrics and hermitian varieties. Secondly by using the projective classification of quadrics in \lbrack 9, p.4\rbrack\  
and their geometric structure, we find some interesting bounds on the number of rational points on the intersection of two quadrics. Thus, we show that for $\mathcal{X}$ a non-degenerate quadric:  $$\# \mathcal{X}_{ Z(\mathcal{Q})}(\mathbb{F}_{q} ) \le 2q^2+3q+1$$ is the best possible. We also show that for $\mathcal{X}$ degenerate and rank $(\mathcal{X})=3$, we get $$\# \mathcal{X}_{ Z(\mathcal{Q})}(\mathbb{F}_{q} ) \le 4q^2+q+1,$$ and this bound is optimal. Identically, for rank$(\mathcal{X})=4$ and g$(\mathcal{X})$=2, we have the optimal bound: $$\# \mathcal{X}_{ Z(\mathcal{Q})}(\mathbb{F}_{q} ) \le 4q^2+1.$$ For $\mathcal{X}$  of rank$(\mathcal{X})=4$ and g$(\mathcal{X})$=1, we get the bound  $$\# \mathcal{X}_{ Z(\mathcal{Q})}(\mathbb{F}_{q} ) \le 3q^2+q+1.$$ Next, we express the exact parameters of the codes $C_2(\mathcal{X})$, the geometric structure of the minimum weight codewords and  show that the performances of the codes $C_2(\mathcal{X})$ defined on the non-degenerate quadrics are better than the ones defined on the degenerate quadrics.\\
In the section 5, by using again the projective classification of quadrics in \lbrack 9, p.4\rbrack , those of hermitian surfaces in \lbrack 8, p.112\rbrack, some properties of hermitian varieties, we find the first five best upper bounds:
$$\# \mathcal{X}_{ Z(\mathcal{Q})}(\mathbb{F}_{q} )=2t^5+t^3+2t^2+1, 2t^5+t^3+t^2+1, 2t^5+2t^2+1,$$$2t^5+t^2+1, 2t^5-t^3+2t^2+1,$\\ on the number of rational points on the intersection of the non-singular hermitian variety $\mathcal{X}$ and any quadric $\mathcal{Q}$.  In section 6, when $\mathcal{X}$ is a non-singular hermitian variety, we determine exactly the minimum distance of the code $C_2(\mathcal{X})$ and the minimum weight codewords. We also list the second, third, fourth, fifth weight, and describe the corresponding codewords reaching these weights.\\
The paper ends with two conjectures: the first, on the minimum distance for the code $C_h(\mathcal{X})$ where $\mathcal{X} \subset {\mathbb{P}}^4(\mathbb{F}_q)$  is the non-singular hermitian variety and the second, on the distribution of the codewords of the first five weights of $C_2(\mathcal{X})$ when $\mathcal{X}\subset\mathbb{P}^N(\mathbb{F}_q)$ is the non-singular hermitian variety.
 \section{Generalities}
We denote by $\mathbb{F}_q$ the field with $q$ elements. Let $V=A^{n+1}(\mathbb{F}_q)$ be the affine space of dimension $n+1$ over $\mathbb{F}_q$ and ${\mathbb{P}}^{n}(\mathbb{F}_q)=\Pi_n$ the corresponding projective space. Then $$\pi_n=\#{\mathbb{P}^{n}(\mathbb{F}_q)}=\#\Pi_n=q^n+q^{n-1}+...+1.$$
We use the term forms of degree two to describe homogeneous polynomials $f$ of degree two, and $\mathcal{Q}=Z(f)$ (the zeros of $f$ in the projective space $\mathbb{P}^{n}(\mathbb{F}_q)$) is a quadric.\\
Let $\mathcal{Q}$ be a quadric. The rank of $\mathcal{Q}$, denoted r$(\mathcal{Q})$, is the smallest number of indeterminates appearing in $f$ under any change of coordinate system. The quadric $\mathcal{Q}$ is said to be degenerate if r$(\mathcal{Q})<n+1$; otherwise it is non-degenerate. For $\mathcal{Q}$ a degenerate quadric and r$(\mathcal{Q})$=r, $\mathcal{Q}$ is a cone $\Pi_{n-r}\mathcal{Q}_{r-1}$ with vertex $\Pi_{n-r}$ (the set of singular points of $\mathcal{Q}$) and base $\mathcal{Q}_{r-1}$ in a subspace $\Pi_{r-1}$ skew to $\Pi_{n-r}$. For degenerate hermitian varieties we have an analogous decomposition as for quadrics. 
\begin{definition}\label{quadrique type} For $\mathcal{Q}=\Pi_{n-r}\mathcal{Q}_{r-1}$ a degenerate quadric with r$(\mathcal{Q})=r$,  $\mathcal{Q}_{r-1}$ is called the non-degenerate quadric associated to $\mathcal{Q}$.
The degenerate quadric $\mathcal{Q}$ will said to be of type hyperbolic (resp. elliptic, parabolic) if its associated non-degenerate quadric is of the same type.
\end{definition}
\begin{definition}\lbrack 11, p.158\rbrack\   
Let $\mathcal{Q}_{1}$ and $\mathcal{Q}_{2}$ be two quadrics. The order $w(\mathcal{Q}_{1}, \mathcal{Q}_{2})$ of the pair $(\mathcal{Q}_{1}$, $\mathcal{Q}_{2})$, is the minimum number of variables, after invertible linear change of variables, necessary to write $\mathcal{Q}_{1}$ and $\mathcal{Q}_{2}$.
\end{definition}
\begin{definition}
For any variety $\mathcal{V}$, the maximum dimension g$(\mathcal{V}$) of linear subspaces lying on $\mathcal{V}$, is called the projective index of $\mathcal{V}$.
\end{definition}
First of all, let us recall the classification of quadrics in $\mathbb{P}^4(\mathbb{F}_q)$. It can be found in the book of J. W. P. Hirschfeld \lbrack 9, p.4\rbrack.
  \hspace{0.5mm}
 \begin{table}[htdp]
 \begin{center}
\begin{tabular}{|c|c|c|c|}

 	\hline
	 r($\mathcal{Q}$) & Description & $\vert \mathcal{Q}\vert $ &g$(\mathcal{Q})$ \\
  	\hline
	\hline
	1 		      & repeated hyperplane $\Pi_3\mathcal{P}_0$ 			  &$q^3+q^2+q+1$&3				 							       						  \\
	\hline
	2  & 	pair of distinct hyperplanes $\Pi_{2}\mathcal{H}_1$ 	  & $2q^3+q^2+q+1$ &3				 							       					\\
	
	\hline		
	2 		      & the plane ${\Pi}_{2}\mathcal{E}_1$	  			  & $q^2+q+1$&2			 							       		 \\
	\hline
	3	      & the cone  ${\Pi}_{1}\mathcal{P}_2$ 			  & $(q+1)(q^2+1)$&2				 							       			\\
	\hline
	 	4	       & the cone ${\Pi}_{0}\mathcal{H}_3(\mathcal{R}, \mathcal{R}^{\prime})$			  & $q(q+1)^2+1$&2				 							       				 \\
	\hline

	 4             &    the cone ${\Pi}_{0}\mathcal{E}_3$                               &$q(q^2+1)+1$ &1                                                        
	
	                                    \\
	
     \hline  
        5      &the parabolic quadric $\mathcal{P}_4$ &     $( q+1)(q^2+1)$&1
                                                           \\
   \hline          
   \end{tabular}
	\end{center}
	\caption{ Quadrics in $\mathbb{P}^{4}(\mathbb{F}_{q})$}
\end{table}%


From the work of Ray-Chaudhuri \lbrack 13, pp.132-136\rbrack\, we deduce that a non-degenerate quadric $\mathcal{P}_4$ in $ \mathbb{P}^4(\mathbb{F}_q)$ contains exactly $\alpha_q=\pi_3$ lines and there are exactly $q+1$ lines contained in $\mathcal{P}_4$ through a given point in $\mathcal{P}_4$.\\
A regulus defined in the work of Hirschfeld \lbrack 8, p.4\rbrack, is the set of transversals of three skew lines in $\mathbb{P}^3(\mathbb{F}_q)$; it consists of $q+1$ skew lines. Thus, a hyperbolic quadric $\mathcal{Q}$ in $\mathbb{P}^3(\mathbb{F}_q)$ is a pair of complementary reguli (each one of the two reguli generates the whole hyperbolic quadric). It is denoted by $\mathcal{Q}={\mathcal{H}_{3}}(\mathcal{R}, \mathcal{R}^{\prime})$ where $\mathcal{R}$ and $\mathcal{R}^{\prime}$ are the two reguli. \\
Let $\mathcal{F}_{h}$ be the vector space of forms of degree $h$ in $V=\mathbb{A}^{n+1}(\mathbb{F}_q)$, $X \subset \mathbb{P}^{n}(\mathbb{F}_q)$ a variety and $\vert X\vert$ the number of rational points of $X$ over $\mathbb{F}_q$. We denote by $W_i$ the set of points with homogeneous coordinates $(x_0:...:x_n)\in \mathbb{P}^n(\mathbb{F}_q)$ such that $x_j=0$ for $j<i$ and $x_i \ne 0$. The family $\{W_i\}_{0\le i\le n}$ is a partition of $\mathbb{P}^n( \mathbb{F}_q)$.
The code $C_h(X)$ is the image of the linear map  
 $c: \mathcal{F}_{h}
  \longrightarrow
  \mathbb{F}_{q}^{\vert X\vert}$, defined by $c(f)={(c_x(f)})_{x\in X}$,  where  $c_x(f)= f(x_0,...,x_n)/{x_i}^h$ with $x=(x_0:...:x_n) \in W_i$. The length of $C_h(X)$ is equal to $\#X(\mathbb{F}_q)$. The dimension of $C_h(X)$ is equal to $\dim{{\mathcal{F}}_h}-\dim{\ker c}$. Therefore, when $c$ is injective we get: 
  \begin{equation}
  \label{dimducode}
  \dim\ {C_{h}(\mathcal{X})} =\left( 
  \begin{array}{c}
  n+h\\
  h
\end{array}
\right).
\end{equation}
The minimum distance of $C_h(X)$ is equal to the mimimum over all $f$ of $\#X(\mathbb{F}_q)-\#{X}_{ Z(f)}(\mathbb{F}_{q})$. 	
\section{Intersection of two quadrics}
In this section, we will estimate the number of points in the intersection of two quadrics $\mathcal{X}$ and $\mathcal{Q}$ in $\mathbb{P}^4(\mathbb{F}_q)$. Some bounds in $\mathbb{P}^4(\mathbb{F}_q)$ have been given by Y. Aubry \lbrack 1\rbrack, G. Lachaud  \lbrack 10\rbrack, D. B. Leep and L. M. Schueller \lbrack 11\rbrack. These bounds are not the best possible even in the case of $\mathbb{P}^4(\mathbb{F}_q)$. Here $\mathcal{X}$ denotes a quadric  of rank$(\mathcal{X})=3,4,5$ in $\mathbb{P}^4(\mathbb{F}_q)$.\\
Let us recall the result of D. B. Leep and L. M. Schueller \lbrack 11, p.172\rbrack\  which gives an upper bound for the number of intersection of a pair of quadrics. 
\begin{theorem}
Let $\mathcal{Q}_1$ and $\mathcal{Q}_2$ be two quadrics in $\mathbb{P}^n(\mathbb{F}_q)$.
Suppose that $w(\mathcal{Q}_1,\mathcal{Q}_2)=n+1$. If $n+1\ge 5$ and $n+1$ odd, then 
$$\vert \mathcal{Q}_1\cap \mathcal{Q}_2\vert \le (2q^{n-1}-q^{n-2}+q^{\frac{n+2}{2}}-q^{\frac{n}{2}}-1)/(q-1).$$ 
\end{theorem}
From this theorem, we deduce the following result.
\begin{corollary} \label{Leep} Let $\mathcal{Q}_1$ and $\mathcal{Q}_2$ be two quadrics in $\mathbb{P}^4(\mathbb{F}_q)$ with $w(\mathcal{Q}_1,\mathcal{Q}_2)=5$. Then, $$\vert \mathcal{Q} _1\cap \mathcal{Q}_2\vert \le 3q^2+q+1.$$
\end{corollary} 
In the theorem above the hypothesis on the order of the two quadrics $\mathcal{Q}_1$ and $\mathcal{Q}_2$ is too restrictive and does not work in general case. Let us state another property, free from this condition. 
\begin{lemma}\label{Serre quadrique}
Let $\mathcal{Q}_1$ and $\mathcal{Q}_2$ be two quadrics in $\mathbb{P}^n(\mathbb{F}_q)$ and $l$ an integer such that $1\le l\le n-1$.\\
Suppose that $w(\mathcal{Q}_1,\mathcal{Q}_2)=n-l+1$ (i.e. there exists a linear transformation such that $\mathcal{Q}_1$ and $\mathcal{Q}_2$ are defined with the indeterminates $x_0, x_1,...,x_{n-l}$) and $\vert \mathcal{Q}_1\cap \mathcal{Q}_2\cap{\mathbb{E}}_{n-l}(\mathbb{F}_q)\vert \le m$ where $\mathbb{E}_{n-l}(\mathbb{F}_q)=\{x\in \mathbb{P}^n(\mathbb{F}_q)\ \vert x_{n-l+1}=...=x_n=0\}$.
Then $$\vert \mathcal{Q}_1\cap \mathcal{Q}_2\vert \le mq^{l}+\pi_{l-1}.$$
This bound is optimal as soon as $m$ is optimal in $\mathbb{E}_{n-l}(\mathbb{F}_q)$.
\end{lemma}
\textbf{Proof}
If $w(\mathcal{Q}_1,\mathcal{Q}_2)=n-l+1$, then there exists a linear transformation such that $\mathcal{Q}_1$ and $\mathcal{Q}_2$ are defined with the indeterminates $x_0, x_1,...,x_{n-l}$.
Let $m_0\le m$ such that $\vert \mathcal{Q}_1\cap \mathcal{Q}_2\cap{\mathbb{E}}_{n-l}(\mathbb{F}_q)\vert = m_0$.
From the fact that $\mathcal{Q}_1\cap \mathcal{Q}_2$ has exactly $m_0$ projective zeros in  $\mathbb{E}_{n-l}(\mathbb{F}_q)$, we deduce that it has exactly $m(q-1)+1$ affine zeros in $\tilde{\mathbb{E}}(\mathbb{F}_q)$ (the corresponding affine space of dimension $n-l+1$).\\
Let $\tilde{y}=(y_0,y_1,...,y_{n-l})$ be a zero of $\mathcal{Q}_1\cap\mathcal{Q}_2$ in $\tilde{\mathbb{E}}(\mathbb{F}_q)$. We get that $y=(y_0,y_1,...,y_{n-l},x_{n-l+1},...,x_{n})$ is a zero of $\mathcal{Q}_1\cap \mathcal{Q}_2$ in $\mathbb{A}^{n+1}(\mathbb{F}_q)$ for any $(x_{n-l+1},...,x_{n}) \in \mathbb{F}^{l}_q$. And reciprocally every zero of $\mathcal{Q}_1\cap\mathcal{Q}_2$ in $\mathbb{A}^{n+1}(\mathbb{F}_q)$ comes from a zero of $\mathcal{Q}_1\cap \mathcal{Q}_2$ in $\tilde{\mathbb{E}}(\mathbb{F}_q)$. Therefore, we deduce that $\mathcal{Q}_1\cap \mathcal{Q}_2$ contains exactly $\{m_0(q-1)+1\}q^l$ affine zeros in $\mathbb{A}^{n+1}(\mathbb{F}_q)$. Thus, in  $\mathbb{P}^{n}(\mathbb{F}_q)$, we get $\vert \mathcal{Q}_1\cap \mathcal{Q}_2\vert = (\{m_0(q-1)+1\}q^l-1)/(q-1)=m_0q^l+\pi_{l-1}$. 
\begin{lemma}\label{ordre rang}
For $\mathcal{Q}_1$ and $\mathcal{Q}_2$ two quadrics in $\mathbb{P}^n(\mathbb{F}_q)$, we get   
$$w(\mathcal{Q}_1,\mathcal{Q}_2)\ge \sup\{\mathrm{r}(\mathcal{Q}_1), \mathrm{r}(\mathcal{Q}_2)\}.$$
\end{lemma}
\textbf{Proof} We have $\mathrm{r}(\mathcal{Q}_{i})\le w(\mathcal{Q}_1,\mathcal{Q}_2)$ for $i=1, 2$.
\subsection{Plane section of the quadric $\mathcal{X}$: g($\mathcal{Q}$)=2}
Here we study the section of the quadric $\mathcal{X}$ by a quadric $\mathcal{Q}$ with  g($\mathcal{Q}$)=2. We have three cases: $\mathcal{Q}={\Pi}_{2}\mathcal{E}_1$, $\mathcal{Q}= {\Pi}_{1}\mathcal{P}_2$  and $\mathcal{Q}={\Pi}_{0}\mathcal{H}_3(\mathcal{R}, \mathcal{R}^{\prime})$. Here $\mathcal{Q}= {\Pi}_{1}\mathcal{P}_2$ and $\mathcal{Q}={\Pi}_{0}\mathcal{H}_3(\mathcal{R}, \mathcal{R}^{\prime})$ are respectively a set of $q+1$ planes through a common line or a common point.    
\subsubsection{$\mathcal{Q}$ is a quadric with  r($\mathcal{Q}$)=2: $\mathcal{Q}={\Pi}_{2}\mathcal{E}_1$}  Here $\mathcal{Q}$ consists of one plane of $\mathbb{P}^4(\mathbb{F}_q)$. If this plane is not contained in $\mathcal{X}$, then $\mathcal{Q}\cap \mathcal{X}$ is a plane quadric. In the book of J. W. P. Hirschfeld \lbrack 7, p.156\rbrack, we have the following classsification of plane quadrics.

\hspace{0.5mm}
\begin{table}[htdp]
\begin{center}

\hspace{20mm}

\begin{tabular}{|c|c|c|c|}
\hline
	r($ \mathcal{Q}^{\prime}$) & Description & $\vert \mathcal{Q}^{\prime} \vert $ &g($\mathcal{Q}^{\prime}$) \\
  	\hline
	\hline
	1 		      & repeated line $\Pi_1\mathcal{P}_0$  			  &$q+1$		&1		 							       						  \\
	\hline
	2  & pair of lines  $\Pi_0\mathcal{H}_1$ 	  &  $2q+1$	&1				 							       					\\
	
	\hline		
	2 		      & point $\Pi_0\mathcal{E}_1$ 			  & 1	&0			 							       		 \\
	\hline
	3	      & parabolic 	$\mathcal{P}_2$		  & $q+1$	& 0			 							       			\\
	\hline
\end{tabular}
\end{center}
	\caption{Plane quadrics in $\mathbb{P}^{2}(\mathbb{F}_q)$}
\end{table}%

Therefore we deduce that  if $\mathcal{Q}$ is a plane of $\mathbb{P}^4(\mathbb{F}_q)$:\\
\textendash\   If it is contained in $\mathcal{X}$, we get  $\vert \mathcal{X}\cap \mathcal{Q}\vert =q^2+q+1$.\\
\textendash\   Otherwise, we get from the table above that  $\vert \mathcal{X} \cap \mathcal{Q}\vert \le 2q+1$.
\subsubsection{$\mathcal{Q}$ is a quadric with r($\mathcal{Q}$)=3,4: $\mathcal{Q}= {\Pi}_{1}\mathcal{P}_2$ and $\mathcal{Q}={\Pi}_{0}\mathcal{H}_3$}
For $\mathcal{X}$ a non-degenerate quadric, there is no plane in $\mathcal{X}$. Thus each of the $q+1$ planes of  $\mathcal{Q}$ intersects $\mathcal{X}$ in at most $2q+1$ points. Therefore, we deduce that $\vert \mathcal{X} \cap \mathcal{Q}\vert\le 2q^2+3q+1$.\\\\ 
For $\mathcal{X}$ degenerate of r$(\mathcal{X})=3,4$, let us consider the order $w(\mathcal{X}, \mathcal{Q})$ of the pair $\{\mathcal{X}, \mathcal{Q}\}$. We have three possibilities $w(\mathcal{X}, \mathcal{Q})=3,4,5:$
\begin{itemize}
\item[{(i)}]If $w(\mathcal{X}, \mathcal{Q})=5$, from the corollary \ref{Leep} we get $\vert \mathcal{X} \cap \mathcal{Q}\vert \le 3q^2+q+1$.
\item[{(ii)}] If $w(\mathcal{X}, \mathcal{Q})=4$, then there exists a linear transformation such that $\mathcal{X}$ and $\mathcal{Q}$ are defined with the indeterminates $x_0, x_1, x_2, x_3$. \\ Suppose rank($\mathcal{X}$)=3. From \lbrack 5, IV-D-3,4 \rbrack, we get $\vert\mathcal{X}\cap \mathcal{Q}\cap\mathbb{E}_3(\mathbb{F}_q)\vert=4q+1$ or $\vert\mathcal{X}\cap \mathcal{Q}\cap\mathbb{E}_3(\mathbb{F}_q)\vert\le 3q$. From lemma \ref{Serre quadrique} we deduce that either $\vert \mathcal{X}\cap \mathcal{Q}\vert= 4q^2+q+1$  in the case where $\mathcal{X}$ and $\mathcal{Q}$ have exactly four common planes through a line or $\vert \mathcal{X}\cap \mathcal{Q}\vert\le 3q^2+1$ otherwise.\\  If rank($\mathcal{X}$)=4 and $g(\mathcal{X})=1$, from \lbrack 5, IV-D-2\rbrack \ we get  $\vert\mathcal{X}\cap \mathcal{Q}\cap\mathbb{E}_3(\mathbb{F}_q)\vert\le 2(q+1)$. Therefore lemma \ref{Serre quadrique} gives $\vert \mathcal{X}\cap \mathcal{Q}\vert\le 2q^2+2q+1$.\\
If rank($\mathcal{X}$)=4 and $g(\mathcal{X})=2$, \lbrack 5, IV-D-3 and IV-E-1\rbrack \  give $\vert\mathcal{X}\cap \mathcal{Q}\cap\mathbb{E}_3(\mathbb{F}_q)\vert=4q$ or  $\vert\mathcal{X}\cap \mathcal{Q}\cap\mathbb{E}_3(\mathbb{F}_q)\vert\le 3q+1$. In fact, in the case $\vert\mathcal{X}\cap \mathcal{Q}\cap\mathbb{E}_3(\mathbb{F}_q)\vert=4q$, $\mathcal{X}\cap \mathcal{Q}\cap\mathbb{E}_3(\mathbb{F}_q)$ is exactly a set of four lines with two lines in each regulus of $\mathcal{H}_3$. Therefore, from lemma \ref{Serre quadrique} and the geometric structure of $\mathcal{X}$ and $\mathcal{Q}$, we get $\vert \mathcal{X}\cap \mathcal{Q}\vert= 4q^2+1$ when $\mathcal{X}$ and $\mathcal{Q}$ have the same vertex $\Pi_{0}$ and contain four common planes with the following configuration:\\
(the first two planes ${\Pi_2}^{(1)}$ and ${\Pi_2}^{(2)}$ meet at the point $\Pi_{0}$, and the two others ${\Pi_2}^{(3)}$ and ${\Pi_2}^{(4)}$ also meet at the point $\Pi_{0}$,\\
the plane ${\Pi_2}^{(1)}$ meets ${\Pi_2}^{(3)}$ and ${\Pi_2}^{(4)}$ respectively at two distinct lines $\mathcal{D}_{1,3}$ and $\mathcal{D}_{1,4}$,\\
the plane ${\Pi_2}^{(2)}$ meets ${\Pi_2}^{(3)}$ and ${\Pi_2}^{(4)}$ respectively at two distinct lines $\mathcal{D}_{2,3}$ and $\mathcal{D}_{2,4}$,\\
each one of the four lines $\mathcal{D}_{1,3}$, $\mathcal{D}_{1,4}$, $\mathcal{D}_{2,3}$ and $\mathcal{D}_{2,4}$ passes through the point $\Pi_{0}$).
For example, for the two quadrics defined by $f_{\mathcal{X}}=x_0x_1+x_2x_3$ and $f_{\mathcal{Q}}=x_3x_0+x_1x_2$ we have $\vert \mathcal{X}\cap \mathcal{Q}\vert= 4q^2+1$. Otherwise, we get $\vert \mathcal{X}\cap \mathcal{Q}\vert \le 3q^2+q+1$.      
\item[{(iii)}] If $w(\mathcal{X}, \mathcal{Q})=3$, necessary from lemma \ref{ordre rang}, we get r($\mathcal{X}$)=r($\mathcal{Q}$)=3. Here $\mathcal{Q}\cap\mathbb{E}_2(\mathbb{F}_q)$ and $\mathcal{X}\cap\mathbb{E}_2(\mathbb{F}_q)$ are two irreducible curves (conics). Therefore from  the theorem of B\'ezout (the fact that two plane conics have exactly four common points or less than three points) and lemma \ref{Serre quadrique}, we deduce that $\vert \mathcal{X}\cap \mathcal{Q}\vert= 4q^2+q+1$, or otherwise $\vert \mathcal{X}\cap \mathcal{Q}\vert\le 3q^2+q+1$.
\end{itemize}
\subsection{Hyperplane section of the quadric $\mathcal{X}$: g ($\mathcal{Q}$)=3}
This paragraph deals with the section of $\mathcal{X}$ by $\mathcal{Q}$ in the case g($\mathcal{Q}$)=3 i.e. $\mathcal{Q}$ contains a hyperplane. We have two cases: r($\mathcal{Q}$)=1 (i.e. $\mathcal{Q}$ is a repeated hyperplane), or r($\mathcal{Q}$)=2 with $\mathcal{Q}$ a pair of hyperplanes.
\subsubsection{$\mathcal{Q}$ is a quadric with r$(\mathcal{Q})$=1: $\mathcal{Q}=\Pi_3\mathcal{P}_0$}
Here $\mathcal{Q}$ is a repeated hyperplane $H$ .
Let us recall three general results on hyperplane section of quadrics in  $\mathbb{P}^n(\mathbb{F}_q)$. 
\begin{lemma}\label{Swin}Swinnerton-Dyer \lbrack 18, p.264\rbrack,      
Let  $\tilde{\mathcal{X}}\subset \mathbb{P}^n(\mathbb{F}_q)$ be a degenerate quadric of rank $r$ and $H$ an hyperplane. Then $\tilde{\mathcal{X}}\cap H$ is a quadric in $\mathbb{P}^{n-1}(\mathbb{F}_q)$ of rank $r, r-1$, or $r-2$.
\end{lemma}
\begin {lemma}\label{Prim} Primrose \lbrack12, pp.299-300\rbrack,
 Let $\tilde{\mathcal{X}}$ be a non-degenerate quadric in $\mathbb{P}^n(\mathbb{F}_q)$  and $H$ a hyperplane.
If $H$ is tangent to $\tilde{\mathcal{X}}$, then $\tilde{\mathcal{X}}\cap H $ is a degenerate quadric of rank $n-1$ in  $\mathbb{P}^{n-1}(\mathbb{F}_q)$. If $H$ is not tangent to $\tilde{\mathcal{X}}$, then $\tilde{\mathcal{X}}\cap H $ is a non-degenerate quadric in  $ \mathbb{P}^{n-1}(\mathbb{F}_q)$. 
\end{lemma}
Let $H \subset \mathbb{P}^4(\mathbb{F}_{q})$ be a hyperplane. If $\mathcal{X}$ is a non-degenerate quadric, from lemma \ref{Prim}, we get 
\begin{equation*}
(\star)\quad \# \mathcal{X}_{H}(\mathbb{F}_{q})=
         \begin{cases}
          {(q+1)}^2, q^2+1&    \text{\  \ if\  H\  is\  not\  tangent \ to \  $\mathcal{X}$,} \\
            q^2+q+1 &  \text{ \ if \ H \ is \ tangent \ to \  $\mathcal{X}$}.
            \end{cases} 
\end{equation*}
If $\mathcal{X}$ is a degenerate quadric, from lemma \ref{Swin}, we get 
\begin{equation}
\label{3.2.1.1}
 \# \mathcal{X}_{H}(\mathbb{F}_{q})\le 2q^2+q+1.
 \end{equation}
The third important result on hyperplane section of a quadric is: 
\begin {lemma}\label{Wolf} J. Wolfmann \lbrack19, pp.191-192\rbrack,
Let $\tilde{\mathcal{X}}$ be a non-degenerate quadric in $\mathbb{P}^n(\mathbb{F}_q)$  and $H$ a hyperplane.
If $H$ is tangent to $\tilde{\mathcal{X}}$, then $\tilde{\mathcal{X}}\cap H$ is of the same type as $\tilde{\mathcal{X}}$ (in the sense of the definition \ref{quadrique type}).
\end{lemma}
\subsubsection{$\mathcal{Q}$ is a quadric with r$(\mathcal{Q})$=2 and $\mathcal{Q}$ is a pair of hyperplanes}
Let $H_1$ and $H_2$ be the two distinct hyperplanes generating $\mathcal{Q}$. We have $\mathcal{Q}=H_1\cup H_2$. Let $\mathcal{P}= H_1\cap H_2$ be the plane of intersection of the two hyperplanes. 
We have 
\begin{equation}
\label{3.2.2.1}
\vert \mathcal{Q }\cap \mathcal{X} \vert = \vert  H_{1} \cap \mathcal{X} \vert  +\vert H_{2} \cap \mathcal{X} \vert - \vert \mathcal {P}  \cap \mathcal{X} \vert.
\end{equation}
Let $\hat{\mathcal{X}_1}=H_1\cap \mathcal{X}$, and $\hat{\mathcal{X}_2}=H_2\cap \mathcal{X}$. We have 
\begin{equation}
\label{3.2.2.2}
\mathcal{P}\cap \mathcal{X}= \mathcal{P}\cap \hat{\mathcal{X}_1}= \mathcal{P}\cap \hat{\mathcal{X}_2}.
\end{equation}
\paragraph{ 3.2.2.1 If $\mathcal{X}$ is non-degenerate (i.e. parabolic)} 
\begin{itemize}
\item[{(i)}]In the case where each hyperplane is tangent to $\mathcal{X}$, we know that $\hat{\mathcal{X}_1}$, and $\hat{\mathcal{X}_2}$ are quadric cones of rank 3. We get that $\vert \hat{\mathcal{X}_1}\vert=\vert \hat{\mathcal{X}_2}\vert=q^2+q+1$. From lemma \ref{Swin} and table 2, we deduce that $\vert\mathcal{P}\cap \hat{\mathcal{X}_1}\vert \ge 1$. Therefore, from relation \ref{3.2.2.1}, we deduce that $\vert \mathcal{X}\cap \mathcal{Q}\vert \le  2q^2+2q+1$.
\item[{(ii)}]In the case where one hyperplane $H_2$ is tangent to $\mathcal{X}$, and the second hyperplane $H_1$ is non-tangent to $\mathcal{X}$, $\hat{\mathcal{X}_1}$ is a non-degenerate quadric surface in $\mathbb{P}^3(\mathbb{F}_{q})$: elliptic or hyperbolic. \\  
\textendash If $\hat{\mathcal{X}_1}$ is an elliptic quadric, from lemmas \ref{Prim}, \ref{Wolf} and table 2, we get that $\mathcal{P}\cap \hat{\mathcal{X}_1}$ is either a plane conic (parabolic), or a single point according to $\mathcal{P}$ being non-tangent or tangent to  $\hat{\mathcal{X}_1}$. Therefore, from relations (\ref{3.2.2.2}), (\ref{3.2.2.1}) and table 2, we deduce that $\vert \mathcal{X}\cap \mathcal{Q}\vert \le  2q^2+q+1$.\\
\textendash If $\hat{\mathcal{X}_1}$ is an hyperbolic quadric, from lemmas \ref{Prim}, \ref{Wolf} and table 2, we get that $\mathcal{P}\cap \hat{\mathcal{X}_1}$ is either a plane conic (parabolic), or a pair of two distinct lines according to $\mathcal{P}$ being non-tangent or tangent to  $\hat{\mathcal{X}_1}$. Therefore, from (\ref{3.2.2.2}) and (\ref{3.2.2.1}), we deduce that $\vert \mathcal{X}\cap \mathcal{Q}\vert \le  2q^2+q+1$.
\item[{(iii)}]In the case where each hyperplane is non-tangent to $\mathcal{X}$, we know that $\hat{\mathcal{X}_1}$, and $\hat{\mathcal{X}_2}$ are non-singular quadric surfaces: elliptic or hyperbolic. They can be of the same type or of different types. \\
\textendash If one of the two quadrics is elliptic, from lemma \ref{Prim} and table 2, we deduce that $\vert\mathcal{P}\cap \hat{\mathcal{X}_1}\vert \ge1$. Therefore from (\ref{3.2.2.1}) and $(\star)$ we get that $\vert \mathcal{X}\cap \mathcal{Q}\vert \le  2q^2+2q+1$.\\  
\textendash If the two quadrics are hyperbolic, from lemmas \ref{Prim}, \ref{Wolf} and table 2, we deduce that $\vert\mathcal{P}\cap \hat{\mathcal{X}_1}\vert \ge q+1$. Therefore, from (\ref{3.2.2.2}) and (\ref{3.2.2.1}), we get $\vert \mathcal{X}\cap \mathcal{Q}\vert\le 2q^2+3q+1$. And this upper bound is reached when $\mathcal{P}$ is non-tangent to $\hat{\mathcal{X}_1}$ and $\hat{\mathcal{X}_2}$. For example, for the two quadrics defined by $f_{\mathcal{X}}=x_0x_1+x_2x_3+x_{4}^{2}$ and $f_{\mathcal{Q}}=(x_0+x_1)(x_2+x_3)$ with $\mathrm{car}(\mathbb{F}_q)\ne 2$ 
\end{itemize}
\paragraph{3.2.2.2 If $\mathcal{X}$ is degenerate (i.e. rank ($\mathcal{X}$)=3,4)}.\\
 Here we have $w(\mathcal{X}, \mathcal{Q})=3,4,5:$
\begin{itemize}
\item[{(i)}] If $w(\mathcal{X}, \mathcal{Q})=5$, from corollary \ref{Leep} we get $\vert \mathcal{X} \cap \mathcal{Q}\vert \le 3q^2+q+1$.
\item[{(ii)}] If $w(\mathcal{X}, \mathcal{Q})=4$. 
Suppose rank($\mathcal{X}$)=3. From  \lbrack 5, IV-C \rbrack \ we get $\vert\mathcal{X}\cap \mathcal{Q}\cap\mathbb{E}_3(\mathbb{F}_q)\vert=4q+1$ or $\vert\mathcal{X}\cap \mathcal{Q}\cap\mathbb{E}_3(\mathbb{F}_q)\vert\le 3q+1$. From lemma \ref{Serre quadrique} we deduce that either $\vert \mathcal{X}\cap \mathcal{Q}\vert= 4q^2+q+1$ in the case where $\mathcal{Q}$ is union of two hyperplanes (non-tangent) each through a pair of planes of $\mathcal{X}$ and the plane of intersection of the two hyperplanes intersecting $\mathcal{X}$ in a line or $\vert \mathcal{X}\cap \mathcal{Q}\vert\le 3q^2+q+1$ otherwise.\\ If rank ($\mathcal{X}$)=4 and $g(\mathcal{X})=1$ from \lbrack 5, IV-B\rbrack \  we get $\vert\mathcal{X}\cap \mathcal{Q}\cap\mathbb{E}_3(\mathbb{F}_q)\vert\le 2(q+1)$. Therefore lemma \ref{Serre quadrique} gives  $\vert \mathcal{X}\cap \mathcal{Q}\vert\le 2q^2+2q+1$.\\
If rank ($\mathcal{X}$)=4 and $g(\mathcal{X})=2$, the result of \lbrack 5, IV-B \rbrack \  gives $\vert\mathcal{X}\cap \mathcal{Q}\cap\mathbb{E}_3(\mathbb{F}_q)\vert=4q$ or $\vert\mathcal{X}\cap \mathcal{Q}\cap\mathbb{E}_3(\mathbb{F}_q)\vert\le 3q+1$. Therefore from lemma \ref{Serre quadrique}, we get  either $\vert \mathcal{X}\cap \mathcal{Q}\vert= 4q^2+1$ when each hyperplane is tangent to $\mathcal{X}$ with the plane of intersection meeting $\mathcal{X}$ at two lines or $\vert \mathcal{X}\cap \mathcal{Q}\vert \le 3q^2+q+1$ otherwise.
\item[{(iii)}] If $w(\mathcal{X}, \mathcal{Q})=3$, necessary from lemma \ref{ordre rang},
 we get r($\mathcal{X}$)=3. Here, the quadric $\mathcal{Q}$ describes a pair of lines in $\mathbb{E}_2(\mathbb{F}_q)$. The number of points in the intersection of two secant lines with a conic (non-singular plane quadric) is exactly four or less than three. From table 2 and lemma \ref{Serre quadrique} we get that either $\vert \mathcal{X}\cap \mathcal{Q}\vert = 4q^2+q+1$ or $\vert \mathcal{X}\cap \mathcal{Q}\vert\le 3q^2+q+1$.
 \end{itemize}
\subsection{ Line section of the quadric $\mathcal{X}$:  g($\mathcal{Q}$)=1}
In this section we estimate the number of points in the intersection of $\mathcal{X}$ and a quadric $\mathcal{Q}$ with g$(\mathcal{Q})=1$. In this case $\mathcal{Q}$ is the degenerate quadric $\Pi_{0}\mathcal{E}_3$ or the non-degenerate quadric (i.e. parabolic) $\mathcal{P}_4$.
\subsubsection{$\mathcal{Q}$ is a quadric with r$(\mathcal{Q})$=4: $\mathcal{Q}=\Pi_{0}\mathcal{E}_3$}
\begin{itemize}
\item[{(i)}]For $\mathcal{X}$ non-degenerate quadric, we get two possibilities:\\
\textendash \ If there is a line of $\mathcal{X}\cap \mathcal{Q}$ through the vertex $\Pi_0$ of the cone $\mathcal{Q}=\Pi_{0}\mathcal{E}_3$, it is obvious that this vertex is a point of $\mathcal{X}$. Therefore there are at most $q+1$ lines of the cone $\mathcal{Q}$ through $\Pi_{0}$ contained in $\mathcal{X}$; the other lines of $\mathcal{Q}$ meet $\mathcal{X}$ in at most two points each. Thus, 
we get $\vert \mathcal{X}\cap \mathcal{Q}\vert\le 2q^2+2q+1$.\\
\textendash\  If there is no line of $\mathcal{X}\cap \mathcal{Q}$ through the vertex of the cone $\mathcal{Q}=\Pi_{0}\mathcal{E}_3$, each line of $\mathcal{Q}$ intersecting $\mathcal{X}$ in at most two points, we deduce that  $\vert \mathcal{X}\cap \mathcal{Q}\vert\le 2(q^2+1)$.
\item[{(ii)}]For $\mathcal{X}$ degenerate quadric (i.e. r$(\mathcal{X})=3,4$).\\
From lemma \ref{ordre rang}, we have $w(\mathcal{X}, \mathcal{Q})=4,5$. The case $w(\mathcal{X}, \mathcal{Q})=5$ follows from corollary \ref{Leep}. Let us consider now that $w(\mathcal{X}, \mathcal{Q})=4$.\\
If r$(\mathcal{X})$=3,4 and g$(\mathcal{X})$=2, or r$(\mathcal{X})$=4 and g$(\mathcal{X})$=1, from \lbrack5, IV-D-2, IV-E-2\rbrack\ we get $\vert\mathcal{X}\cap \mathcal{Q}\cap\mathbb{E}_3(\mathbb{F}_q)\vert\le 2(q+1)$. Thus, we deduce from lemma \ref{Serre quadrique} that $\vert \mathcal{X}\cap \mathcal{Q}\vert\le 2q^2+2q+1$.
\end{itemize}
\begin{Remark}Let $\mathcal{Q}$ be a non-degenerate quadric (i.e. $\mathcal{Q}=\mathcal{P}_4$) and r$(\mathcal{X})$=3,4.
The cases r$(\mathcal{X})$=3,4 and g$(\mathcal{X})$=2, correspond to the first part of 3.1.2. with $\mathcal{Q}$ at the place of $\mathcal{X}$. In the same way r$(\mathcal{X})$=4 and g$(\mathcal{X})$=1, corresponds to 3.3.1.(i).
\end{Remark}
\subsubsection{Intersection of two non-degenerate quadrics}
Here we study the number of points in the intersection of two non-degenerate quadrics $\mathcal{X}$ and $\mathcal{Q}$.
\begin{Proposition}\label{fred}
Let $\mathcal{X}$ and $\mathcal{Q}$ be two non-degerate quadrics in $\mathbb{P}^4(\mathbb{F}_q)$. If there is no line in $\mathcal{X}\cap \mathcal{Q}$, then  $\vert \mathcal{X}\cap\mathcal{Q}\vert \le 2(q^2+1)$.     
\end{Proposition}
\textbf{Proof} We use the fact that $\alpha_q$ is the number of lines of $\mathcal{X}$, each one of them meeting $\mathcal{Q}$ in at most two points. Moreover there pass exacly $q+1$ lines through each point contained in $\mathcal{X}$.\\\\
Now we will study the section of two non-degenerate quadrics containning a common line. 
Let us consider a line $\mathcal{D}$ contained in $\mathcal{X}\cap \mathcal{Q}$ and $\mathcal{P}$ a plane through $\mathcal{D}$. By the principle of duality in the projective space \lbrack15, pp.49-51\rbrack, or \lbrack7, p.33, theorem 3.1 p.85\rbrack, we deduce that there are exactly $q+1$ hyperplanes $(H_i)$ through  $\mathcal{P}$. These $q+1$ hyperplanes $(H_i)$  generate $\mathbb{P}^4(\mathbb{F}_q)$.\\
For $i=1,..., q+1$, we denote $\hat{\mathcal{X}}_i=H_i \cap \mathcal{X}$ and   $\hat{\mathcal{Q}}_i=H_i \cap \mathcal{Q}$.
Thus, we get: 
\begin{equation}
\label{uneseuledroite}
 \vert \mathcal{X} \cap \mathcal{Q} \vert \le \vert\hat{\mathcal{X}}_1\cap \hat{\mathcal{Q}}_1\vert + \sum_{i=2}^{q+1}\vert (\hat{\mathcal{X}}_i \cap \hat{\mathcal{Q}}_i)-\mathcal{D}\vert.
\end{equation} 
From lemma \ref{Prim}, we deduce that $\hat{\mathcal{X}}_i$ and $\hat{\mathcal{Q}}_i$ are quadrics of rank 3 or 4 in $\mathbb{P}^3(\mathbb{F}_q)$. Since they contain the line $(\mathcal{D})$, they can not be elliptic.  They are either  hyperbolic or cone quadrics (of rank 3). Thus one has to study the three types of intersection in $\mathbb{P}^3(\mathbb{F}_q)$ given by the table 3.\\

\hspace{0.5mm}
\begin{table}[htdp]
\begin{center}

\hspace{15mm}
 
\begin{tabular}{|c|c|c|}

	\hline
	Types &  $\hat{\mathcal{X}_i}\cap \hat{\mathcal{Q}_i}$  \\
  	\hline
	\hline
	1 		      & $(\mathrm{hyperbolic\  quadric}) \cap (\mathrm{quadric\ cone})$			  				 							       						  \\
	\hline
            2             &   $(\mathrm{quadric\  cone}) \cap (\mathrm{quadric\ cone})$					 							       					\\

	\hline		
	3 		     &  $(\mathrm{hyperbolic\  quadric}) \cap (\mathrm{hyperbolic\  quadric})$			 							       		 \\
	\hline
		 		\end{tabular}
\end{center}
	\caption{Intersection of $\hat{\mathcal{X}_i}\cap \hat{\mathcal{Q}_i}$  in $\mathbb{P}^3(\mathbb{F}_q)$}
\end{table}%
\begin{lemma}\label{deux hyperplans}
Let $\mathcal{X}$ and $\mathcal{Q}$  be two non-degerate quadrics in $\mathbb{P}^n(\mathbb{F}_q)$, and $\mathcal{K}$ a linear space of codimension 2.
Then there exists at most two hyperplanes $H_i$, $i=1,2$ through  $\mathcal{K}$ such that $\hat{\mathcal{X}_i}= \ \hat{\mathcal{Q}_i}$. 
\end{lemma}
\textbf{Proof} Let $$\mathcal{Q}=\sum_{0\le i\le j\le n}^{}{a_{ij}}{x_i}{x_j}{}\ \mathrm{and} \ \mathcal{X}=\sum_{0\le i\le j\le n}^{}{a_{ij}^{\prime}}{x_i}{x_j}{}.$$ 
Without loss of generality, we can choose a system of coordinates, such that $H_1=\{x_n=0\}$, $H_2=\{x_{n-1}=0\}$. 
For $H_1\cap \mathcal{X}=  H_1\cap \mathcal{Q}$, we get that $a_{ij}= a_{ij}^{\prime}$ except may be for $(i,n)$ with $i=0,...,n$. 
For $H_2\cap \mathcal{X}=  H_2\cap \mathcal{Q}$, we get that $a_{ij}= a_{ij}^{\prime}$ except maybe for $(i,n-1)$ with $i=0,...,n-1$ and $(n-1,n)$. Therefore, we deduce that $a_{ij}= a_{ij}^{\prime}$ except for $(i,j)=(n-1,n)$; and there exists $(\alpha,\beta)\in \mathbb{F}_{q}^{2}$ ($\alpha \ne \beta$) such that:

\begin{equation*}
   {}
         \begin{cases}
         \mathcal{Q}= \mathcal{Q}_{0}(x_0, x_1,..., x_{n-1}, x_{n})+ \alpha x_{n-1}x_{n} &\text {}\\
        \mathcal{X}= \mathcal{Q}_{0}(x_0, x_1,..., x_{n-1},x_{n})+ \beta x_{n-1}x_{n}. &\text {}
\end{cases} 
\end{equation*} 
Let us suppose that there is a third hyperplane $H_{3}$ through $\mathcal{K}$ such that $H_3\cap \mathcal{X}= H_3\cap \mathcal{Q}$. We can suppose that $H_3=\{ax_{n-1}+bx_{n}=0\}$ and $b\ne 0$, therefore $x_n=-\frac{a}{b}x_{n-1}$. From the above system and the fact that $H_3\cap\mathcal{X}=H_3\cap \mathcal{Q}$ we deduce that $-\alpha\frac{a}{b}x_{n-1}^2=-\beta\frac{a}{b}x_{n-1}^2$. Since $\alpha\ne\beta$, we deduce that $a=0$, which leads to $H_3=H_2$.
\begin{Proposition}
Let $\mathcal{X}$ and $\mathcal{Q}$ be two non-degerate quadrics in $\mathbb{P}^4(\mathbb{F}_q)$ containing a line ($\mathcal{D}$).  
There is a plane ($\mathcal{P}$) containing ($\mathcal{D}$) such that  $\hat{\mathcal{X}_i}\ne \ \hat{\mathcal{Q}_i}$ for $i=1,...,q+1$.
\end{Proposition}
\textbf{Proof}
We know from lemma \ref{deux hyperplans} that there exists at most two hyperplanes $H_i$, $i=1,2$ such that $\hat{\mathcal{X}_i}= \ \hat{\mathcal{Q}_i}$. There are also at most $2(q+1)$ planes through ($\mathcal{D}$) contained in $H_1\cup H_2$. Since there are exactly  $q^2+q+1$ planes through ($\mathcal{D}$), we conclude that there are at least $q^2-q-1>0$ possibilities of choice for the plane ($\mathcal{P}$) such that $\hat{\mathcal{X}_i}\ne \ \hat{\mathcal{Q}_i}$ for $i=1,...,q+1$.
\begin{lemma}\label{unique hyperplan}
Let $\mathcal{D}_1$, $\mathcal{D}_2$ be two secant lines of the non-degenerate quadric $\mathcal{X}\subset \mathbb{P}^4(\mathbb{F}_q)$ and $\mathcal{P}=< \mathcal{D}_1, \mathcal{D}_2>$ the plane defined by these two lines.
Let $(H_i)$ i=1,...,q+1 be the $q+1$ hyperlanes containing $\mathcal{P}$.
If there exists a hyperplane $H_1$ tangent to $\mathcal{X}$, then it is unique (i.e. the remaining $q$ hyperplanes are non-tangent to $\mathcal{X}$). 
\end{lemma}
\textbf{Proof}
Let $H_1$ containing $\mathcal{P}$, be a hyperplane which is tangent to $\mathcal{X}$ at a point $P_1$. Then from lemma \ref{Prim}, $\hat{\mathcal{X}}_1$ is a quadric cone of rank 3 in $\mathbb{P}^3(\mathbb{F}_q)$: $\hat{\mathcal{X}}_1$ is a set of $q+1$ lines through $P_1$ ($\mathcal{D}_1$ and $\mathcal{D}_2$ are two of them).\\ 
Let $H_2$ containing $\mathcal{P}$ be another hyperplane which is tangent to $\mathcal{X}$ at a point $P_2$. Then $\hat{\mathcal{X}}_2$ is a set of $q+1$ lines through $P_2$ and ($\mathcal{D}_1$ and $\mathcal{D}_2$ are two of them). The lines $(\mathcal{D}_1)$ and $(\mathcal{D}_2)$ intersect at $P_1$ and $P_2$, therefore we deduce that $P_1=P_2$. Thus, there exist a unique hyperplane $H_1$ tangent to $\mathcal{X}$ at the point $P=\mathcal{D}_1\cap\mathcal{D}_2$. 
\begin{Remark}\label{quadriques deux droites}
We also know from \lbrack 5,\S IV-D-4 and IV-E-1\rbrack \  that if the intersection of two quadrics cone (of rank 3) or two hyperbolic quadrics in $\mathbb{P}^3(\mathbb{F}_q)$ contains three lines, it contains exactly four common lines.
\end{Remark}
From the results of \lbrack 5, \S IV\rbrack \ and table 3 above, we deduce the following table 4. 
 
\hspace{0.5mm}
\begin{table}[htdp]
\begin{center}
\hspace{15mm}

\begin{tabular}{|c|c|c|c|}

	\hline
	Types & 4 lines & 2 lines  & 1 line\\
  	\hline
	\hline
	1 		      & 			  &$3q$   & $2q+1$				 							       						  \\
	\hline
	2  & 4q+1  	  &  $3q$		& $2q+1$	   		 							       					\\
	
	\hline		
	3 		      & $4q$			  &   $3q+1$  &     $2(q+1)$				 						\\	       		 
	\hline
	\end{tabular}
\end{center}
	\caption{Number of points and lines in $\hat{\mathcal{X}_i}\cap \hat{\mathcal{Q}_i}$           }
\end{table}%
Let us explain the table 4.
For the type 1, the intersection of a hyperbolic quadric and a quadric cone, contains at most two lines; therefore $3q$ and $2q+1$ are respectively the maximum number of this intersection when containing (exactly or) at most  two lines or exactly one line. From the above remark, the types 2 and 3 are describe as before.\\\\
Now an estimation on the number of points in the intersection of the two non-degenerate quadrics  $\mathcal{X}$ and $\mathcal{Q}$ containning a common line is reduced to the two following simple cases.
\paragraph{If $\mathcal{X}\cap\mathcal{Q}$ contains exactly one common line:} Each $\hat{\mathcal{X}_i}\cap \hat{\mathcal{Q}_i}$ contains exacly one line, and from the table 4, we get for $i=1,...,q+1$ $\vert \hat{\mathcal{X}_i}\cap \hat{\mathcal{Q}_i} \vert \le 2(q+1)$. Therefore from relation (\ref{uneseuledroite}), we deduce that $\vert \mathcal{X}\cap \mathcal{Q} \vert \le q^2+3q+2$.
\paragraph{If $\mathcal{X}\cap\mathcal{Q}$ contains at least two common lines:}We distinguish two following cases:
\subparagraph{1. In the case where $\mathcal{X}\cap\mathcal{Q}$ contains only skew lines:} from table 4, we get  for $i=1,...,q+1$ $\vert \hat{\mathcal{X}_i}\cap \hat{\mathcal{Q}_i} \vert \le 2q+1$ for the types 1 and 2. Indeed, if $\hat{\mathcal{X}_i}\cap \hat{\mathcal{Q}_i} $ contains more than one line, from the fact that one of the two quadrics $\hat{\mathcal{X}_i}$ or $ \hat{\mathcal{Q}_i}$ is a quadric cone, they are secant. For type 3, if $\hat{\mathcal{X}_i}\cap \hat{\mathcal{Q}_i}$ contains exactly four common lines, then it contains necessarily two secant lines. This is a contradiction. From remark \ref{quadriques deux droites},   $\hat{\mathcal{X}}_i \cap \hat{\mathcal{Q}}_i$ can only contain at most two lines. Thus, from table 4, $3q+1$ is an upper bound for the number of points in $\hat{\mathcal{X}}_i\cap \hat{\mathcal{Q}}_i$. The lines of $\hat{\mathcal{X}}_i\cap \hat{\mathcal{Q}}_i$ are skew; so they belong to the same regulus. From \lbrack 5, IV-E-1\rbrack \  we get $\vert \hat{\mathcal{X}_i}\cap \hat{\mathcal{Q}_i} \vert \le 2(q+1)$. Finally, from relation (\ref{uneseuledroite}) we deduce that $\vert \mathcal{X} \cap \mathcal{Q} \vert \le q^2+3q+2$.
\subparagraph{2. In the case where there exist some secant lines in $\mathcal{X}\cap\mathcal{Q}$,} let $\{\mathcal{D}_1, \mathcal{D}_2\}$ denote a pair of secant lines and $\mathcal{P}$ be the plane generated by this pair of secants. Let  $$\mathcal{Q}=\sum_{0\le i\le j\le 4}^{}{a_{ij}}{x_i}{x_j}{}\ \mathrm{and} \ \mathcal{X}=\sum_{0\le i\le j\le 4}^{}{a_{ij}^{\prime}}{x_i}{x_j}{}.$$ 
\begin{itemize}         
\item[(i)] If there are exactly two hyperplanes $H_i$ $i=1,2$ such that $\hat{\mathcal{X}}_i= \hat{\mathcal{Q}}_i$, then we can choose a system of coordinates, such that $H_1=\{x_4=0\}$ and $H_2=\{x_3=0\}$.
From the proof of lemma \ref{deux hyperplans}, 
there exists $(\alpha,\beta)\in \mathbb{F}_{q}^{2}$ ($\alpha \ne \beta$) such that:

\begin{equation*}
   {}
         \begin{cases}
         \mathcal{Q}= \mathcal{Q}_{0}(x_0, x_1, x_2, x_3, x_4)+ \alpha x_3x_4 &\text {}\\
        \mathcal{X}= \mathcal{Q}_{0}(x_0, x_1, x_2, x_3,x_4)+ \beta x_3x_4. &\text {}
\end{cases} 
\end{equation*} 
Thus, we have $\mathcal{X}\cap\mathcal{Q}=(H_1\cap \mathcal{X})\cup (H_2\cap \mathcal{X})$ and from the relation (\ref{3.2.2.1}) we get that $\vert \mathcal{X}\cap\mathcal{Q}\vert \le 2q^2+2q+1.$
\item [(ii)] If there is exactly one hyperplane $H_1$ such that  $\hat{\mathcal{X}}_1= \hat{\mathcal{Q}}_1$, then by the same reasoning as above, we can choose two distinct linear forms $h_1(x_0,x_1,x_2, x_3, x_4)$ and $h_2(x_0,x_1,x_2, x_3, x_4)$ such that:  
\begin{equation*}
     {}
         \begin{cases}
         \mathcal{Q}= \mathcal{X}_{0}(x_0, x_1, x_2, x_3)+ x_4h_1(x_1,x_2, x_3,x_4) \text {}\\
         \mathcal{X}= \mathcal{X}_{0}(x_0, x_1, x_2, x_3)+ x_4h_2(x_1,x_2, x_3,x_4). \text {}
\end{cases} 
\end{equation*} 
We have $\mathcal{X}\cap\mathcal{Q}=(H_1\cap \mathcal{X})\cup (H_2\cap\mathcal{X}\cap \mathcal{Q})$  where $H_2$ is the hyperplane defined by the linear form $h_1-h_2$. We also have $\vert H_2\cap\mathcal{X}\cap \mathcal{Q}\vert \le 4q+1$ from table 4 and $\vert H_1\cap\mathcal{X}\vert \le q^2+2q+1$. Therefore we get that $\vert \mathcal{X}\cap\mathcal{Q}\vert \le q^2+6q+2$.
\item [(iii)] If for $i=1,...,q+1$ we get $\hat{\mathcal{X}}_i \ne \hat{\mathcal{Q}}_i$, one has two possibilities: 
If the $q+1$ hyperplanes $H_i$ are all non-tangent to $\mathcal{X}$, $\hat{\mathcal{X}_i}$ are non-degenerate quadrics and only the types 1 and 3 of the table 3 can appear. From table 4  and remark \ref{quadriques deux droites}, we deduce that $\vert \hat{\mathcal{X}_i}\cap \hat{\mathcal{Q}_i} \vert \le 4q$. If there exists a hyperplane $H_i$ tangent to $\mathcal{X}$, then from lemma \ref{unique hyperplan}, it is unique. Let $H_1$ be this tangent hyperplane, $\hat{\mathcal{X}_1}$ is a quadric cone and $\hat{\mathcal{X}_1}\cap \hat{\mathcal{Q}_1}$ is of type 1 or 2. Therefore we get $\vert \hat{\mathcal{X}_1}\cap \hat{\mathcal{Q}_1} \vert \le 4q+1$. Finally from relation
\begin{equation}   
 \label{dedrsecantes}   
   \vert \mathcal{X} \cap \mathcal{Q} \vert \le \vert\hat{\mathcal{X}}_1\cap \hat{\mathcal{Q}}_1\vert + \sum_{i=2}^{q+1}\vert (\hat{\mathcal{X}}_i \cap \hat{\mathcal{Q}}_i)-({\mathcal{D}}_1\cup {\mathcal{D}}_2)\vert
  \end{equation}
we deduce that, $\vert \mathcal{X}\cap \mathcal{Q} \vert \le 2q^2+3q+1$.
 \end{itemize}
\section{The parameters of the code $C_2(\mathcal{X})$ defined on the quadric $\mathcal{X}$}
When $\mathcal{X}=Z(f^{\prime})\subset \mathbb{P}^{n}({\mathbb{F}}_{q})$ is a quadric, the linear map  
 $c: \mathcal{F}_{2}
  \longrightarrow
  \mathbb{F}_{q}^{\vert \mathcal{X}\vert}$ is not injective because $\ker{c}=\{ \lambda f^{\prime} \vert \lambda \in \mathbb{F}_{q} \}$. Therefore we deduce that: $$\dim\ {C_{2}(\mathcal{X})} =\left( 
  \begin{array}{c}
  n+2\\
  2
\end{array}
\right)-1=\frac{n(n+3)}{2}$$
From the results of section 3, we deduce the followings results.
\begin{theorem}
Let $\mathcal{Q}$ be a quadric in $\mathbb{P}^{4}(\mathbb{F}_q)$ and $\mathcal{X}$ the non-degenerate (parabolic) quadric in $\mathbb{P}^{4}(\mathbb{F}_q)$. We get  $$\# \mathcal{X}_{ Z(\mathcal{Q})}(\mathbb{F}_{q} ) \le 2q^2+ 3q+1$$ and this bound is the best possible.\\
The code $C_{2}(\mathcal{X})$ defined on the parabolic quadric $\mathcal{X}$ is a $\lbrack n, k,d \rbrack_{q}$-code where  \\
 $ n= (q+1)(q^2+1)$,  
 $k =14 $,
 $d=q^3-q^2-2q$.
\end{theorem}
\begin{theorem}
The minimum weight codewords of the code $C_2(\mathcal{X})$ correspond to:\\
\textendash either quadrics which are union of two (non-tangent) hyperplanes each intersecting $\mathcal{X}$ at a hyperbolic quadric such that the plane of intersection of the two hyperplanes intersects $\mathcal{X}$ at a plane conic.\\
\textendash or quadrics with r$(\mathcal{Q})=3$ (i.e. $q+1$ planes through a line) and each plane containing exactly two lines of $\mathcal{X}$.\\
\textendash or non-degenerate (i.e. parabolic) quadrics containing two secant lines of $\mathcal{X}$ defining a plane contained in the $q+1$ hyperplanes which has one tangent hyperplane, the $q$ other non-tangent to $\mathcal{X}$, and with maximal hyperplane section.
\end{theorem}
\begin{theorem}
Let $\mathcal{Q}$ be a quadric in $\mathbb{P}^{4}(\mathbb{F}_q)$ and  $\mathcal{X}$ a degenerate quadric in $\mathbb{P}^{4}(\mathbb{F}_q)$.
We get :\\
\textendash If $\mathcal{X}$ is degenerate with r$(\mathcal{X})=3$, then
 $$\# \mathcal{X}_{ Z(\mathcal{Q})}(\mathbb{F}_{q} ) = 4q^2+q+1\quad \mathrm{or} \qquad \# \mathcal{X}_{ Z(\mathcal{Q})}(\mathbb{F}_{q} ) \le 3q^2+q+1.$$ 
\textendash If $\mathcal{X}$ is degenerate with r$(\mathcal{X})$=4 and $g(\mathcal{X})=2$, then
 $$\# \mathcal{X}_{ Z(\mathcal{Q})}(\mathbb{F}_{q} ) = 4q^2+1\quad \mathrm{or} \qquad \# \mathcal{X}_{ Z(\mathcal{Q})}(\mathbb{F}_{q} ) \le 3q^2+q+1.$$ 
   \end{theorem}
\begin{theorem}
The code $C_{2}(\mathcal{X})$ defined on the degenerate quadric $\mathcal{X}$ in PG(4,q) is a $\lbrack n, k,d \rbrack_{q}$-code where:  \\
\textendash If $\mathcal{X}$ is degenerate with r$(\mathcal{X}$)=3, then
 $ n= q^3+q^2+q+1$,  
 $k =14 $,
 $d=q^3-3q^2$.\\
 \textendash If $\mathcal{X}$ is degenerate with r($\mathcal{X}$)=4 and $g(\mathcal{X})=2$, then
 $ n= q^3+2q^2+q+1$,  
 $k =14 $,
 $d=q^3-2q^2+q.$
 \end{theorem}
 \begin{theorem}
The minimum weight codewords of the code $C_2(\mathcal{X})$ correspond :\\
For $\mathcal{X}$ degenerate with r($\mathcal{X}$)=3 to:\\
\textendash quadrics which are union of two hyperplanes (non-tangent) each through a pair of planes and the plane of intersection of the two hyperplanes intersecting $\mathcal{X}$ in a line.\\
\textendash quadrics with r$(\mathcal{Q})=3$ and a plane containing exactly four lines of $\mathcal{X}$.\\
For $\mathcal{X}$ degenerate quadric with r$(\mathcal{X})=4$ and $g(\mathcal{X})=2$:\\ 
\textendash quadrics which are union of two tangent hyperplanes to $\mathcal{X}$ and the plane of intersection of the two hyperplanes meeting $\mathcal{X}$ at two secant lines.\\
\textendash quadrics with r$(\mathcal{Q})=4$ and $g(\mathcal{Q})=2$  containning exactly four common planes of $\mathcal{X}$ with the following configuration: (the two first planes ${\Pi}_2^{(1)}$ and ${\Pi}_2^{(2)}$ meet at the point $\Pi_{0}$, and the two others planes ${\Pi}_{2}^{(3)}$ and ${\Pi}_{2}^{(4)}$ meet at the point $\Pi_{0}$,\\
the plane ${\Pi}_{2}^{(1)}$ meets ${\Pi}_{2}^{(3)}$ and ${\Pi}_{2}^{(4)}$ respectively at two distinct lines $\mathcal{D}_{1,3}$ and $\mathcal{D}_{1,4}$,\\
the plane ${\Pi}_{2}^{(2)}$ meets ${\Pi}_{2}^{(3)}$ and ${\Pi}_{2}^{(4)}$ respectively at two distinct lines $\mathcal{D}_{2,3}$ and $\mathcal{D}_{2,4}$,\\
each one of the four lines $\mathcal{D}_{1,3}$, $\mathcal{D}_{1,4}$, $\mathcal{D}_{2,3}$ and $\mathcal{D}_{2,4}$ pass through the point $\Pi_{0}$).
\end{theorem}
\begin{theorem}
Let $\mathcal{Q}$ be a quadric in $\mathbb{P}^{4}(\mathbb{F}_q)$ and  $\mathcal{X}$ a degenerate quadric in $\mathbb{P}^{4}(\mathbb{F}_q)$ of rank $(\mathcal{X})=4$ with $g(\mathcal{X})=1$.
We get :
 $$\# \mathcal{X}_{ Z(\mathcal{Q})}(\mathbb{F}_{q} ) \le 3q^2+q+1.$$ 
 The code $C_{2}(\mathcal{X})$ defined on the degenerate quadric $\mathcal{X}$ of rank($\mathcal{X}$)=4 with $g(\mathcal{X})=1$ is a $\lbrack n, k,d \rbrack_{q}$-code where 
 $ n= q^3+q+1$,  
 $k =14 $,
 $d \ge q^3-3q^2$.
\end{theorem}
\begin{Remark}
From the study of the parameters of these codes, we can assert that the performances of the codes $C_2(\mathcal{X})$ defined on the non-degenerate quadrics are better than the ones defined on the degenerate quadrics.
\end{Remark}
\section{Quadratic section of the non-degenerate hermitian variety}
In this section $\mathbb{F}_q$ denotes the field with $q$ elements, where $q=t^2$ and $\mathcal{X}$ denotes the non-degenerate (i.e. non-singular) hermitian variety of $\mathbb{P}^4(\mathbb{F}_q) $ of equation $\mathcal{X}: x_0^{t+1}+x_1^{t+1}+x_2^{t+1}+x_3^{t+1}+x_4^{t+1}=0$.\\
In $\lbrack 3,\mathrm{p}.1175\rbrack$ Bose and Chakravarti proved the following result:
\begin{theorem}\label{cardinal hermitienne}
Let $\tilde{\mathcal{X}}\subset \mathbb{P}^n(\mathbb{F}_{q})$ be a non-degenerate hermitian variety. Then, 
\begin{equation}
\label{bochanbrehermitian}
\#\tilde{\mathcal{X}}(\mathbb{F}_{q})=\Phi(n,t^2)=\lbrack t^{n+1}-(-1)^{n+1}\rbrack\lbrack t^{n}-(-1)^{n}\rbrack/(t^2-1)
\end{equation}
For $\tilde{\mathcal{X}}\subset \mathbb{P}^n(\mathbb{F}_{q})$ a degenerate hermitian variety of rank $r<n+1$, we have:
 $$\#\tilde{\mathcal{X}}(\mathbb{F}_{q})=(t^{2}-1)\pi_{n-r}\Phi(r-1,t^2)+\pi_{n-r}+\Phi(r-1,t^2),$$
where $\Phi(n,t^2)$ is given by (\ref{bochanbrehermitian}).
\end{theorem}
From theorem \ref{cardinal hermitienne} we get that $$\#\mathcal{X}(\mathbb{F}_q)=(t^2+1)(t^5+1).$$ It has been shown by Bose et al. $\lbrack 3,\mathrm{p}. 1176\rbrack$ that g($\mathcal{X})=1$ (i.e. $\mathcal{X}$ contains lines and does not contain a plane). In the work of Hirschfeld $\lbrack 9,\mathrm{p}.60 \rbrack$ we get that there are $t^3+1$ lines contained in $\mathcal{X}$ through each point of $\mathcal{X}$.\\ 

We will study here the section of $\mathcal{X}$ by any quadric $\mathcal{Q}$ of  $\mathbb{P}^4(\mathbb{F}_q)$. An approach to this problem has also been considered by F. Rodier in \lbrack 14, pp.207-208\rbrack. He used the result of G. Lachaud \lbrack 10, proposition 2.3\rbrack \  which gives an upper bound for the number of rational points of an algebraic set of a given degree. He found a bound which is not optimal.\\
We will use a different method depending mainly on the geometric structure of the quadric $\mathcal{Q}$ and the hermitian variety $\mathcal{X}$.\\
Let us recall the classification of hermitian varieties in $\mathbb{P}^3(\mathbb{F}_q)$, table.5, which can be found in the work of J.W.P. Hirschfeld $\lbrack 9,\mathrm{p}.60 \rbrack$. \\
 \hspace{0.5mm}
 \begin{table}[htdp]
 \begin{center}
\hspace{15mm}

\begin{tabular}{|c|c|c|c|}

	\hline
	r($\mathcal{V}$) & Description & $\vert\mathcal{V}\vert $ &g($\mathcal{V}$) \\
  	\hline
	\hline
	1 		      & repeated plane $\Pi_{2}\mathcal{U}_0$ 			  &$t^4+t^2+1$	&2			 							       						  \\
	\hline
	2  & $t+1$ collinear planes  	$\Pi_{1}\mathcal{U}_1$	  &  $t^5+t^4+t^2+1$		&2			 							       					\\
	
	\hline		
	3 		      & a cone 	$\Pi_{0}\mathcal{U}_2$		  & $t^5+t^2+1$	&1			 							       		 \\
	\hline
	4	      & non-singular hermitian surface 	$\mathcal{U}_3$		  & $t^5+t^3+t^2+1$	&1			 							       			\\
	\hline
\end{tabular}
\end{center}
	\caption{Hermitian surfaces in $\mathbb{P}^{3}(\mathbb{F}_q)$}
\end{table}
\subsection{Plane section of the non-degenerate hermitian variety $\mathcal{X}$: g($\mathcal{Q}$)=2}
We are interested in the case  where g($\mathcal{Q}$)=2 for the quadric $\mathcal{Q}$. These are the cases where $\mathcal{Q}={\Pi}_{2}\mathcal{E}_1$, $\mathcal{Q}= {\Pi}_{1}\mathcal{P}_2$  or $\mathcal{Q}={\Pi}_{0}\mathcal{H}_3(\mathcal{R}, \mathcal{R}^{\prime})$. In the book of J. W. P. Hirschfeld \lbrack 7, p.160\rbrack, we have also the classification of plane hermitian curves. 
\hspace{0.5mm}
 \begin{table}[htdp]
 \begin{center}

\hspace{15mm}

\begin{tabular}{|c|c|c|c|}

	\hline
	r($\mathcal{V}$) & Description & $\vert\mathcal{V}\vert $ &g($\mathcal{V}$) \\
  	\hline
	\hline
	1 		      & repeated line  $\Pi_{1}\mathcal{U}_0$			  &$t^2+1$	&1			 							       						  \\
	\hline
	2  & cone	$\Pi_{0}\mathcal{U}_1$  	  &  $t^3+t^2+1$	&1				 							       					\\
	
	\hline		
	3 		      & non-singular hermitian curve $\mathcal{U}_2$			  & $t^3+1$	&0			 							       		 \\
	\hline
		 		\end{tabular}
\end{center}
	\caption{Plane hermitian curves}
\end{table}%
\\
From the table above, for the section of the non-singular hermitian variety $\mathcal{X}$ by a plane $\mathcal{Q}$ we get $\vert \mathcal{X}\cap \mathcal{Q}\vert \le t^3+t^2+1$. In the case  r($\mathcal{Q})=3$ or r($\mathcal{Q})=4$, ($\mathcal{Q}$ is a union of $q+1$ planes) we get $\vert \mathcal{X}\cap \mathcal{Q}\vert \le t^5+t^4+t^3+2t^2+1$.
\subsection{Hyperplane section of the non-degenerate hermitian variety $\mathcal{X}$: g($ \mathcal{Q}$)=3}
In the case g($\mathcal{Q}$)=3, i.e. when $\mathcal{Q}$ contains a hyperplane, we have two possibilities: r($\mathcal{Q}$)=1 ($\mathcal{Q}$ is a hyperplane), or r($\mathcal{Q}$)=2 and $\mathcal{Q}$ is a pair of distinct hyperplanes.
\subsubsection{$\mathcal{Q}$ is a quadric with r$(\mathcal{Q})$=1: $\mathcal{Q}=\Pi_3\mathcal{P}_0$}
Here $\mathcal{Q}$ is a repeated hyperplane $H$. Let us recall two general results of Bose and Chakravarti on hyperplane section of a non-degenerate hermitian variety.
\begin{theorem} \lbrack 3, p.1173\rbrack\ \label{BC1}
Let $\tilde{\mathcal{X}}$ be a non-degenerate hermitian variety in $\mathbb{P}^n(\mathbb{F}_q)$  and $H$ a hyperplane.
If $H$ is tangent to $\tilde{\mathcal{X}}$ at $P$, then $\tilde{\mathcal{X}}\cap H $ is a degenerate hermitian variety  of rank $n-1$ in  $\mathbb{P}^{n-1}(\mathbb{F}_q)$. The singular space of $\tilde{\mathcal{X}}\cap H$ consists of the single point $P$.
\end{theorem}
\begin{theorem}\lbrack 4, p.272\rbrack\   \label{Char}
Let $\tilde{\mathcal{X}}$ be a non-degenerate hermitian variety in $\mathbb{P}^n(\mathbb{F}_q)$  and $H$ a hyperplane.
If $H$ is not tangent to $\tilde{\mathcal{X}}$, then $\tilde{\mathcal{X}}\cap H $ is a non-degenerate variety in  $\mathbb{P}^{n-1}(\mathbb{F}_q)$. 
\end{theorem}
From theorems \ref{cardinal hermitienne}, \ref{BC1} and \ref{Char}, we deduce the following result. 
\begin{theorem}\label{BoseChar P4}
Let $H \subset \mathbb{P}^4(\mathbb{F}_{q})$ be a hyperplane
\begin{equation*}
 \# \mathcal{X}_{H}(\mathbb{F}_{q})=
         \begin{cases}
          t^5+ t^3+t^2+1 &    \text{\  \ if\  H\  is\  not\  tangent \ to \  $\mathcal{X}$,} \\
            t^5+t^2+1 &  \text{ \ if \ H \ is \ tangent \ to \  $\mathcal{X}$}.
            \end{cases} 
\end{equation*}                         
\end{theorem}
If $H$ is tangent to $\mathcal{X}$, then $\mathcal{X}\cap H$ is a singular hermitian surface of rank 3 in $\mathbb{P}^3(\mathbb{F}_{q})$. More precisely, $\mathcal{X}\cap H$ is a set of $t^3+1$ lines passing through a common point $P$. If $H$ is non-tangent to $\mathcal{X}$, then $\mathcal{X}\cap H$ is a non-singular hermitian surface in $\mathbb{P}^3(\mathbb{F}_{q})$; thus we get $\#  \mathcal{X}_{H}(\mathbb{F}_{q})=t^5+t^3+t^2+1$.
\subsubsection{$\mathcal{Q}$ is a quadric with r$(\mathcal{Q})$=2 and $\mathcal{Q}$ is a pair of hyperplanes}
Here we will also use two important results of Bose and Charkravarti on degenerate hermitian varieties in $\mathbb{P}^{n}(\mathbb{F}_q)$.
\begin{theorem} \lbrack 3, p.1171\rbrack\ \label{BC2}
Let $\tilde{\mathcal{X}}$ be a degenerate hermitian variety of rank $r<n+1$ in $\mathbb{P}^n(\mathbb{F}_q)$  and $\Pi_{r-1}$ a linear projective space of dimension $r-1$ disjoint from the singular space $\Pi_{n-r}$ of $\tilde{\mathcal{X}}$.
Then $\Pi_{r-1}\cap \tilde{\mathcal{X}}$ is a non-degenerate hermitian variety in $\Pi_{r-1}$. 
\end{theorem}
\begin{theorem} \lbrack 3, p.1171\rbrack\ \label{BC3}
Let $\tilde{\mathcal{X}}\subset \mathbb{P}^n(\mathbb{F}_{q})$ be a degenerate hermitian variety of rank $r<n+1$. If $P$ is any point belonging to the singular space of $\tilde{\mathcal{X}}$  and $D$ is an arbitrary point of $\tilde{\mathcal{X}}$, then any point of the line $(PD)$ belongs to  $\tilde{\mathcal{X}}$.  
\end{theorem}
Now let $H_1$ and $H_2$ be two distinct hyperplanes generating $\mathcal{Q}$. We have $\mathcal{Q}=H_1\cup H_2$. Let $\mathcal{P}= H_1\cap H_2$ be the plane of intersection of the two hyperplanes. Let $\hat{\mathcal{X}_1}=H_1\cap \mathcal{X}$, and $\hat{\mathcal{X}_2}=H_2\cap \mathcal{X}$.
\begin{itemize}
\item[{(i)}] In the case where each hyperplane is tangent to $\mathcal{X}$,
 we know that $\hat{\mathcal{X}_1}$ and $\hat{\mathcal{X}_2}$ as singular hermitian surfaces are sets of $t^3+1$ lines passing respectively through the point $P_1$ and $P_2$ ($P_1\ne P_2$).
From theorem \ref{BC1}, we deduce that $P_1$ and $P_2$ are the singular spaces of $\hat{\mathcal{X}_1}$ and $\hat{\mathcal{X}_2}$.\\ 
\textendash If the plane $\mathcal{P}$ does not pass through at least one of the two points $P_1$ and $P_2$, without loss of generality we can suppose that $\mathcal{P}$ does not pass through $P_1$. Since $P_1$ is the singular space of $\hat{\mathcal{X}_1}$, from theorem \ref{BC2} we deduce that $\mathcal{P}\cap \hat{\mathcal{X}_1}$ is a non-singular curve in $\mathbb{P}^2(\mathbb{F}_{q})$. Therefore, we get $ \vert\mathcal{P}\cap \hat{\mathcal{X}_1}\vert=t^3+1$, which with the relations (\ref{3.2.2.1}) and (\ref{3.2.2.2}) give $\vert \mathcal{X }\cap \mathcal{Q} \vert=2t^5-t^3+2t^2+1$.\\ 
\textendash If the plane $\mathcal{P}$ passes through $P_1$ and $P_2$, then it is obvious that 
$P_1$ and $P_2$ belong to $\mathcal{P}\cap\hat{\mathcal{X}_1}$ and $\mathcal{P}\cap\hat{\mathcal{X}_2}$ respectively. From the relation (\ref{3.2.2.2}), we deduce that $P_2\in\mathcal{P}\cap\hat{\mathcal{X}_1}$. Therefore $P_2\in\hat{\mathcal{X}_1}$. Since $P_1$ is the singular space of $\hat{\mathcal{X}_1}$, from theorem \ref{BC3} we deduce that the line ($P_1P_2$) is contained in $\hat{\mathcal{X}_1}$. The line ($P_1P_2$) is also contained in $\mathcal{P}$. Therefore ($P_1P_2$) is contained in $\mathcal{P}\cap\hat{\mathcal{X}_1}$ and from table 6, we get that $\mathcal{P}\cap \hat{\mathcal{X}_1}$ is a degenerate hermitian curve in $\mathbb{P}^2(\mathbb{F}_{q})$. Here $\mathcal{P}\cap \hat{\mathcal{X}_1}$ cannot be a set of $t+1$ lines through a common point, because it would imply $P_1=P_2$. Thus, $\mathcal{P}\cap \hat{\mathcal{X}_1}$ is the line defined by the two points $P_1$ and $P_2$. Therefore, from the relations $(\ref{3.2.2.1})$ and (\ref{3.2.2.2}) we get $\vert \mathcal{X }\cap \mathcal{Q}\vert = 2t^5+t^2+1$.
\item[{(ii)}] In the case where one hyperplane ($H_2$) is tangent to $\mathcal{X}$ and the second hyperplane ($H_1$) is non-tangent to $\mathcal{X}$, $\hat{\mathcal{X}_2}$ and   $\hat{\mathcal{X}_1}$ are singular and non-singular hermitian surfaces respectively. From theorems \ref{BC1}, \ref{BC2} and relation (\ref{3.2.2.2}) we deduce that the following two conditions are equivalent:\\
\textendash the plane $\mathcal{P}$ is non-tangent to $\hat{\mathcal{X}_1}$\\
\textendash the plane $\mathcal{P}$ is disjoint from the singular space $\{P_2\}$ of  $\hat{\mathcal{X}_2}$\\ 
Therefore, from the theorems above and relation (\ref{3.2.2.1}) we get either $\vert \mathcal{X }\cap \mathcal{Q}\vert  = 2t^5+2t^2+1$ when $\mathcal{P}$ is non-tangent to $\hat{\mathcal{X}_1}$ or $\vert \mathcal{X }\cap \mathcal{Q}\vert  = 2t^5+t^2+1$  when $\mathcal{P}$ is tangent to $\hat{\mathcal{X}_1}$.
  \item[{(iii)}] In the case where both hyperplanes are non-tangent to $\mathcal{X}$, we know that $\hat{\mathcal{X}_1}$ and $\hat{\mathcal{X}_2}$ are both non-singular hermitian surfaces. From the relation (\ref{3.2.2.2}) and the fact that $\mathcal{P}\cap \hat{\mathcal{X}_1}$ (resp. $\mathcal{P}\cap \hat{\mathcal{X}_2}$) is singular if and only if $\mathcal{P}$ is tangent to $\hat{\mathcal{X}_1}$ (resp. $\hat{\mathcal{X}_2}$), we deduce that $\mathcal{P}$ is either tangent to $\hat{\mathcal{X}_1}$ and $\hat{\mathcal{X}_2}$, or non-tangent to $\hat{\mathcal{X}_1}$ and $\hat{\mathcal{X}_2}$.\\
\textendash If the plane $\mathcal{P}$ is tangent to $\hat{\mathcal{X}_1}$ (i.e. it is also tangent to $\hat{\mathcal{X}_2}$), then from theorem \ref{BC1}, $\mathcal{P}\cap \hat{\mathcal{X}_1}$ is a singular hermitian curve of rank 2 in $\mathbb{P}^2(\mathbb{F}_{q})$. Therefore, we get $\vert\mathcal{P}\cap \hat{\mathcal{X}_1}\vert=t^3+t^2+1$, which with the relation (\ref{3.2.2.1}) give $\vert \mathcal{X }\cap \mathcal{Q} \vert=2t^5+t^3+t^2+1$ \\
 \textendash If the plane $\mathcal{P}$ is non-tangent to $\hat{\mathcal{X}_1}$ (i.e. it is also non-tangent to $\hat{\mathcal{X}_2}$), $\mathcal{P}\cap \hat{\mathcal{X}_1}$ is a non-singular hermitian curve in $\mathbb{P}^2(\mathbb{F}_{q})$. Therefore, we get $ \vert\mathcal{P}\cap \hat{\mathcal{X}_1}\vert=t^3+1$, which with the relation (\ref{3.2.2.1}) give $\vert \mathcal{X }\cap \mathcal{Q} \vert=2t^5+t^3+2t^2+1$. 
\end{itemize}
\subsection{Line section of the non-degenerate hermitian variety $\mathcal{X}$: g($\mathcal{Q}$)=1} 
We will now study the section of the non-degenerate hermitian variety $\mathcal{X}$ by $\mathcal{Q}$ with g($\mathcal{Q}$)=1. In this case $\mathcal{Q}$ is the non-degenerate quadric (parabolic) or the degenerate quadric ${\Pi}_{0}\mathcal{E}_3$.\\         
For the section of $\mathcal{X}$ by $\mathcal{Q}$ we need to distinguish two cases.
\subsubsection{$\mathcal{Q}\cap \mathcal{X}$ does not contain any line}
When $\mathcal{Q}$ is degenerate (i.e. $\mathcal{Q}={\Pi}_{0}\mathcal{E}_3$) and $\mathcal{Q}\cap \mathcal{X}$ does not contain any line, every line of $\mathcal{Q}$ intersects $\mathcal{X}$ in at most $t+1$ points. From the fact that $\mathcal{Q}$ consists of $q^2+1$ lines through the point $\Pi_{0}$, we deduce that $\vert \mathcal{Q}\cap \mathcal{X}\vert \le t^5+t^4+t+1$.\\
When $\mathcal{Q}$ is non-degenerate, it contains exactly $\alpha_q$ lines and there are $q+1$ lines of $\mathcal{Q}$ through any point of $\mathcal{Q}$. Therefore, we deduce that $\vert \mathcal{Q}\cap \mathcal{X}\vert \le\alpha_q(t+1)/(q+1)$ i.e. $\vert \mathcal{Q}\cap \mathcal{X}\vert \le t^5+t^4+t+1$.
\subsubsection{$\mathcal{Q}\cap \mathcal{X}$ contains a line}
First, we will study the simple case of the section of $\mathcal{X}$ by the degenerate quadric.  Next, we will treat the more technical case, section of $\mathcal{X}$ by the non-degenate quadric.
\paragraph{5.3.2.1\quad When $\mathcal{Q}$ is degenerate:}
We get two possibilities:\\
\textendash\  If there is no line of $\mathcal{X}\cap \mathcal{Q}$ through the vertex $\Pi_{0}$ of the cone ${\Pi}_{0}\mathcal{E}_3$, as in 5.3.1  we deduce that  $\vert \mathcal{X}\cap \mathcal{Q}\vert\le t^5+t^4+t+1$.\\
\textendash \ If there is a line of $\mathcal{X}\cap \mathcal{Q}$ through the vertex $\Pi_{0}$ of the cone, it is obvious that $\Pi_{0}$ is a point of $\mathcal{X}$. Therefore, there are at most $t^3+1$ lines of $\mathcal{Q}$ contained in $\mathcal{X}$. Each one of the other $t^4-t^3$ lines of $\mathcal{Q}$ intersects $\mathcal{X}$ in at most $t+1$ points. Thus, we deduce that $\vert \mathcal{Q}\cap \mathcal{X}\vert \le 2t^5-t^4+t^2+1$.
\paragraph{5.3.2.2\quad When $\mathcal{Q}$ is non-degenerate:}
We will use the same technique as in \S 3.3.2 (intersection of two non-degenerate quadrics). 
Let us consider a line $\mathcal{D}$ contained in $\mathcal{Q}\cap \mathcal{X}$, $\mathcal{P}$ a plane through $\mathcal{D}$ and $(H_i)$ the $q+1$ hyperplanes passing through $\mathcal{P}$ which generate  $\mathbb{P}^4(\mathbb{F}_q)$. As in \S 3.3.2, for $i=1,..., q+1$  $\hat{\mathcal{Q}}_i=H_i \cap \mathcal{Q}$ are either hyperbolic quadrics or cone  quadrics in $\mathbb{P}^3(\mathbb{F}_q)$. From theorems \ref{BC1}, \ref{Char} we get that for $i=1,..., q+1$,  $\hat{\mathcal{X}_i}=H_i\cap \mathcal{X}$ are hermitian surfaces of rank 3 or 4 (non-singular hermitian surfaces) in $\mathbb{P}^3(\mathbb{F}_q)$.\\ Thus, one has to study the four types of intersection in $\mathbb{P}^3(\mathbb{F}_q)$ given by the table 7.\\ 
\hspace{0.5mm}
 \begin{table}[htdp]
 \begin{center}
\hspace{15mm}
 
\begin{tabular}{|c|c|c|}

	\hline
	Type &  $\hat{\mathcal{Q}_i}\cap \hat{\mathcal{X}_i}$  \\
  	\hline
	\hline
	1 		      & $(\mathrm{hyperbolic\  quadric}) \cap (\mathrm{non-sing\  herm\ surf})$			  				 							       						  \\
	\hline
            2             &   $(\mathrm{quadric\  cone}) \cap (\mathrm{non-sing \ herm\ surf })$					 							       					\\

	\hline		
	3 		     &  $(\mathrm{hyperbolic\  quadric}) \cap (\mathrm{sing\  herm \ surface})$			 							       		 \\
	\hline
	4	      & $(\mathrm{quadric\  cone}) \cap (\mathrm{sing \ hermitian \ surface})$				 							       			\\
	\hline
	 		\end{tabular}
\end{center}
	\caption{Intersection of $\hat{\mathcal{Q}_i}\cap \hat{\mathcal{X}_i}$ in $\mathbb{P}^{3}(\mathbb{F}_q)$                                         }
\end{table}%

From the paragraphs 4.1 and 4.2 of \lbrack 6\rbrack, we can be more precise about $\# \hat{{\mathcal{X}_i}}_{Z(\hat{\mathcal{Q}_i})}(\mathbb{F}_{q})$  for the types 1 and 2 of table 7 (section of the non-singular hermitian surface $\hat{\mathcal{X}_i}$ by a cone or hyperbolic quadric $\hat{\mathcal{Q}_i}$). Thus, we have the following table 8.\\
\hspace{0.5mm}
 \begin{table}[htdp]
 \begin{center} 
\hspace{5mm}

\begin{tabular}{|c|c|c|l|} 
	\hline
	r$(\hat{\mathcal{Q}}_i)$ & Type& $L(\hat{{\mathcal{X}_i}}\cap \hat{\mathcal{Q}}_i)$   & $\# \hat{{\mathcal{X}_i}}_{Z(\hat{\mathcal{Q}_i})}(\mathbb{F}_{q})$   \\
  	\hline		
	 \hline   
	            
	     3    & 2      &           2         &   $t^3+2t^2-t+1$                                                                                        
                          \\
   \cline{3-4}
  
      (cone)             & {} &        1                    &     $t^3+t^2+1$                                                                    
                            \\

    \hline 
        
    \hline
        {4}          & {}         &            3     &$2t^3+t^2+1$                                                        
                             	   \\          
	           
   \cline{3-4}
      (hyperbolic)       &   1     &      2               &$\le t^3+3t^2-t+1$                                                                                                                                        \\
             \cline{3-4}                              
                       {$\mathcal{H}_{3}( \mathcal{R}_i, \mathcal{R}^{\prime}_i)$}                    &    {}         &   1  &  $\le t^3+2t^2+1$                                       
                                       \\   
          	               \hline

    \end{tabular}
\end{center}
\caption{ Number of points and lines in $\hat{\mathcal{Q}_i}\cap \hat{\mathcal{X}_i}$                  }
\end{table}%

In table 8, $\hat{\mathcal{X}_i}$ and $L(\hat{{\mathcal{X}_i}}\cap \hat{\mathcal{Q}}_i)$ denote respectively the non-degenerate hermitian surface and the number of lines contained in $\hat{{\mathcal{X}_i}}\cap \hat{\mathcal{Q}}_i$. In the case where $\hat{\mathcal{Q}_i}$ is a hyperbolic, we consider $L(\hat{{\mathcal{X}_i}}\cap \mathcal{R}_i)$ (where $\mathcal{R}_i$  is a regulus of the hyperbolic quadric $\hat{\mathcal{Q}_i}$) instead of $L(\hat{{\mathcal{X}_i}}\cap \hat{\mathcal{Q}}_i)$. \\ 

From the tables above, we deduce that an estimation on the number of points in the intersection of $\mathcal{X}\cap\mathcal{Q}$ where $\mathcal{X}$ and $\mathcal{Q}$ are respectively  a non-degenerate hermitian variety and a non-degenerate quadric variety with a common line, is resolved by considering the two following simple cases.
\subparagraph{(a) If $\mathcal{Q}\cap \mathcal{X}$ contains exactly one common line:} From table 8, we get $\vert \hat{\mathcal{Q}_i}\cap \hat{\mathcal{X}_i} \vert \le t^3+2t^2+1$ for the types 1 and 2 of table 7.\\
Likewise, for the type 3, if $\hat{\mathcal{Q}_i}$ and $\hat{\mathcal{X}_i}$ contain exactly the only common line $\mathcal{D}$, a fortiori there is a regulus $\mathcal{R}_{i}$ of the hyperbolic quadric $\hat{\mathcal{Q}_i}$ containing $\mathcal{D}$. Each one of the $q$ other lines of $\mathcal{R}_{i}$ intersects $\hat{\mathcal{X}_i}$ in at most $t+1$ points. We deduce that $\vert \hat{\mathcal{Q}_i}\cap \hat{\mathcal{X}_i} \vert \le t^3+2t^2+1$.\\
For the type 4, the quadric cone $\hat{\mathcal{Q}_i}$ consists of $q+1$ lines through a common point. The fact that there is only one line of $\hat{\mathcal{Q}_i}$ in $\hat{\mathcal{X}_i}$, implies that each one of the $q$ other lines of $\hat{\mathcal{Q}_i}$ intersects $\hat{\mathcal{X}_i}$ in at most $t+1$ points (the vertex of the cone is one of them). We deduce that $\vert \hat{\mathcal{Q}_i}\cap \hat{\mathcal{X}_i} \vert \le t^3+t^2+1$.\\
Finally, when $\mathcal{Q}\cap \mathcal{X}$ contains exactly one common line, we get:
\begin{equation}
 \label{exoncoline}
 \mathrm{for} \quad i=1,2,...,q+1\quad \vert \hat{\mathcal{Q}_i}\cap \hat{\mathcal{X}_i} \vert \le t^3+2t^2+1.
 \end{equation}
In this way, by using the relations (\ref{uneseuledroite}) and (\ref{exoncoline}) we conclude that $\vert \mathcal{Q}\cap \mathcal{X}\vert \le t^5+t^4+t^3+2t^2+1$. 
\subparagraph{(b) If $\mathcal{Q}\cap \mathcal{X}$ contains at least two common lines:} We distinguish two cases: $\mathcal{Q}\cap \mathcal{X}$ containing only skew lines or some secant lines. We will use an important property on the intersection of a hyperbolic quadric and a non-singular hermitian surface in $\mathbb{P}^{3}(\mathbb{F}_q)$.
\begin{lemma}\label{Hirch}\lbrack 8, pp.123-124\rbrack\  
Let $\hat{\mathcal{Q}_{i}}= \mathcal{H}( \mathcal{R}_{i} , \mathcal{R}_{i}^{\prime} )$ and  $\hat{\mathcal{X}_{i}}$ denote respectively the hyperbolic quadric and the non-singular hermitian surface in $\mathbb{P}^{3}(\mathbb{F}_q)$.
If $\hat{\mathcal{Q}_{i}}= \mathcal{H}( \mathcal{R}_{i} , \mathcal{R}_{i}^{\prime} )$ has three skew lines on $\hat{\mathcal{X}_{i}}$, then $\hat{\mathcal{X}_{i}} \cap\hat{\mathcal{Q}_{i}}$ consists of $2(t+1)$ lines of  $\mathcal{R}_{0} \cup  \mathcal{R}^{\prime} _{0} $ where $\mathcal{R}_{0} \subset  \mathcal{R}_{i}$, $\mathcal{R}^{\prime} _{0}  \subset { \mathcal{R}^{\prime}}_{i}$ and  $\vert \mathcal{R}^{\prime} _{0}  \vert = \vert \mathcal{R}_{0} \vert=t+1 $.
\end{lemma}
\begin{Remark}\label{Fred1}
Lemma \ref{Hirch} says that, if a hyperbolic quadric contains three skew lines on the non-singular hermitian surface $\hat{\mathcal{X}_{i}}$, then it contains exactly $2(t+1)$ lines of the surface $\hat{\mathcal{X}_{i}}$, and $t+1$ lines in each of the two reguli.
\end{Remark}
We will also use the following result.
\begin{lemma}\label{unique hyperplan hermitian}
Let $\mathcal{D}_1$, $\mathcal{D}_2$ be two secant lines of the non-degenerate hermitian variety $\mathcal{X}\subset \mathbb{P}^{4}(\mathbb{F}_q)$ and $\mathcal{P}=< \mathcal{D}_1, \mathcal{D}_2>$ the plane defined by these two lines.
Let $(H_i)$ i=1,...,q+1 be the $q+1$ hyperlanes containing $\mathcal{P}$.
If there exists a hyperplane $H_1$ tangent to $\mathcal{X}$, then it is unique (i.e. the remaining $q$ hyperplanes are non-tangent to $\mathcal{X}$). 
\end{lemma}
The proof of this lemma is analogous to the one of lemma \ref{unique hyperplan}.
\begin{itemize}
\item[(i)]
If $\mathcal{Q}\cap \mathcal{X}$ contains only skew lines.\\ 
For the types 2 and 4 of table 7, since the quadric cone consists of $q+1$ lines through a common point, we deduce that $\hat{\mathcal{Q}_i}\cap \hat{\mathcal{X}_i}$ contains exactly one common line. And therefore from table 8, we get that $\vert \hat{\mathcal{Q}_i}\cap \hat{\mathcal{X}_i} \vert \le t^3+t^2+1$ for the type 2. For the type 4, we also have $\vert \hat{\mathcal{Q}_i}\cap \hat{\mathcal{X}_i} \vert \le t^3+t^2+1$ as in the subparagraph (a) above.\\
For the type 3, since the singular hermitian surface $\hat{\mathcal{X}_i}$ is a set of $t^3+1$ lines through a common point, we also deduce that $\hat{\mathcal{Q}_i}\cap \hat{\mathcal{X}_i}$ contains exactly one common line. Thus, as in subparagraph (a), we get $\vert \hat{\mathcal{Q}_i}\cap \hat{\mathcal{X}_i} \vert \le t^3+2t^2+1$.\\
For the type 1, from remark \ref{Fred1}, we deduce that $\hat{\mathcal{Q}_i}$ as a hyperbolic quadric in $\mathbb{P}^{3}(\mathbb{F}_q)$, cannot contain three skew lines of $\hat{\mathcal{X}_i}$. Otherwise $\hat{\mathcal{Q}_i}$ would  contain two secant lines of $\hat{\mathcal{X}_i}$. Thus, under the condition of the non-existence of secant lines in $\mathcal{Q}\cap \mathcal{X}$,  any  regulus of  $\hat{\mathcal{Q}_i}$ contains at most two skew lines of $\hat{\mathcal{X}_i}$ and therefore from the table 8, we get $\vert \hat{\mathcal{Q}_i}\cap \hat{\mathcal{X}_i} \vert \le t^3+3t^2-t+1$.\\
Finally, when $\mathcal{Q}\cap \mathcal{X}$ contains only skew lines, we get: 
\begin{equation}
\label{twcolines}
\mathrm{for} \quad i=1,2,...,q+1\quad \vert \hat{\mathcal{Q}_i}\cap \hat{\mathcal{X}_i} \vert \le t^3+3t^2-t+1.
\end{equation}
In this way, by using the relations (\ref{uneseuledroite}) and (\ref{twcolines}) we conclude that $\vert \mathcal{Q}\cap \mathcal{X}\vert \le t^5+2t^4-3t^2+1$. 
\item[(ii)]
 If $\mathcal{Q}\cap \mathcal{X}$ contains some secant lines, let $(\mathcal{D}_1)$ and $(\mathcal{D}_2)$ be two of them (which are not skew) and $\mathcal{P}$ the plane generated by them. One has two cases:\\
 \textendash \ The $q+1$ hyperplanes $H_i$ are all non-tangent to $\mathcal{X}$: in this case $\hat{\mathcal{X}_i}$ is a non-singular hermitian surface. From the four types of table 7, only the first two types appear in this case. From the table 8, we deduce that 
\begin{equation}
\label{toushyptan}
\mathrm{for} \quad i=1,2,...,q+1\quad \vert \hat{\mathcal{Q}_i}\cap \hat{\mathcal{X}_i} \vert \le 2t^3+t^2+1.
\end{equation}
 Thus, from (\ref{dedrsecantes}) and (\ref{toushyptan}) we get $\vert \mathcal{Q}\cap \mathcal{X}\vert \le 2t^5-t^4+2t^3+t^2+1$. \\
\textendash \  If the $q+1$ hyperplanes are not all non-tangent to $\mathcal{X}$, from lemma \ref{unique hyperplan hermitian}, there is a unique tangent hyperplane to $\mathcal{X}$, say $H_1$. Now, $\hat{\mathcal{Q}_1}\cap \hat{\mathcal{X}_1}$ can be of the types 3 and 4.\\ For the type 3, $\hat{\mathcal{Q}_1}$ and  $\hat{\mathcal{X}_1}$ contain exactly two common lines. Each one of the remaining $t^3-1$ lines of $\hat{\mathcal{X}_1}$ intersects $\hat{\mathcal{Q}_1}$ in at most two points (the vertex of $\hat{\mathcal{X}_1}$ is one of them). We deduce that $\vert \hat{\mathcal{Q}_1}\cap \hat{\mathcal{X}_1} \vert \le t^3+2t^2$.\\ For the type 4, let $t\ge3$, from the B\'ezout theorem, we deduce that deg($\hat{\mathcal{Q}_1}\cap \hat{\mathcal{X}_1})\le 2(t+1)$. Therefore $\hat{\mathcal{Q}_1}\cap \hat{\mathcal{X}_1}$ contains at most $2(t+1)$ lines. Without loss of generality we can assume that $\hat{\mathcal{Q}_1}\cap \hat{\mathcal{X}_1}$ contains exactly $2(t+1)$ lines. If we suppose that  each one of the remaining $t^2-2t-1$ lines of the quadric cone $\hat{\mathcal{Q}_1}$ intersects  $\hat{\mathcal{X}_1}$ at $t+1$ points, we get that  $\vert \hat{\mathcal{Q}_1}\cap \hat{\mathcal{X}_1} \vert \le 3t^3-t+1$. For $t=2$, we get obviously $\vert \hat{\mathcal{Q}_1}\cap \hat{\mathcal{X}_1} \vert \le \vert \hat{\mathcal{Q}_1}\vert = t^4+t^2+1$.\\ Thus, when there is a unique hyperplane tangent to $\mathcal{X}$, from (\ref{dedrsecantes}) we get that: 
\begin{equation*}
 \vert \mathcal{X}\cap{\mathcal{Q}}\vert \le
         \begin{cases}
          2t^5+t^2+1 &    \text{\  \ if\  $t=2$} ,\\
             2t^5-t^4+3t^3-t+1   &  \text{ \ if  \ $t\ge 3$}.
            \end{cases} 
\end{equation*}            
\end{itemize}
\section{The parameters of the code $C_2(\mathcal{X})$ defined on the non-degenerate hermitian variety $\mathcal{X}$}
From the results of section 5, we deduce the following results.
\begin{theorem}
Let $\mathcal{Q}$ be a quadric in $\mathbb{P}^{4}(\mathbb{F}_q)$ and $\mathcal{X}$ the non-degenerate hermitian variety $ \mathcal{X} : x_{0}^{t+1}+x_{1}^{t+1}+x_{2}^{t+1}+x_{3}^{t+1}+x_{4}^{t+1}= 0$, we get  $$\# \mathcal{X}_{ Z(\mathcal{Q})}(\mathbb{F}_{q} ) = 2t^5+ t^{3}+ 2t^{2}+1,\  2t^5+ t^{3}+ t^{2}+1$$   $$ \# \mathcal{X}_{ Z(\mathcal{Q})}(\mathbb{F}_{q} )= 2t^5+ 2t^{2}+1,\  2t^5+ t^{2}+1,$$
\begin{equation*}
\mathrm{and\  for}\ t \neÊ2,3\quad
  \#\mathcal{X}_{ Z(\mathcal{Q})}(\mathbb{F}_{q} )= 
  \begin{cases}
  2t^5-t^3+2t^{2}+1 &  \\ \quad \text{or} \\
   \le 2t^5-t^4+3t^{3}-t+1. & \text{}
   \end{cases}
   \end{equation*}
  \end{theorem}
 \begin{theorem} (Serre-S\o rensen)(\lbrack 16, p.351\rbrack, \lbrack 17, chp.2, pp.7-10\rbrack)
 Let \\$f(x_{0},...,x_{m})$ be a homogeneous polynomial in $m+1$ variables with coefficients in  $\mathbb{F}_{q}$  and degree $h\le q$. Then the number of zeros of $f$ in   $\mathbb{P}^{m}(\mathbb{F}_q)$ satisfies:  $$\#Z_{(f)}(\mathbb{F}_{q})\le hq^{m-1}+\pi_{m-2}.$$
This upper bound is attained when $Z_{(f)}$ is a union of $h$ hyperplanes passing through a common linear space of codimension 2.  
 \end{theorem} 
 \begin{theorem}
The code $C_{2}(\mathcal{X})$ defined on the hermitian variety $\mathcal{X} : x_{0}^{t+1}+x_{1}^{t+1}+x_{2}^{t+1}+x_{3}^{t+1}+x_4^{t+1}= 0$ is a $\lbrack n, k,d \rbrack_{q}$-code where  \\
 $ n= (t^2+1)(t^5+1)$,  
 $k =15 $,
 $d=t^7-t^5-t^3-t^2$.
 \end{theorem}
 \textbf{Proof:}
 Here we need to prove that the linear map  
 $c: \mathcal{F}_{2}
  \longrightarrow
  \mathbb{F}_{q}^{\vert \mathcal{X}\vert}$ is injective. In general, when $f$ is an homogeneous polynomial of degree $h$ in the affine space $\mathbb{A}^{5}(\mathbb{F}_{q})$, by the Serre-S\o rensen's bound on hypersurfaces, we deduce that $\#Z_{(f)}(\mathbb{F}_q)\le ht^6+t^4+t^2+1$. From the fact that $\#\mathcal{X}(\mathbb{F}_q)=t^7+t^5+t^2+1$ we deduce that for $h\le t$, the map $c$ is injective. From the relation (\ref{dimducode}), we get that $k=15$.
From theorem 6.1 and the defintion of $d$ at the end of section 2, we deduce that $d=\#\mathcal{X}(\mathbb{F}_q)-(2t^5+t^3+2t^2+1)=t^7-t^5-t^3-t^2$.\\\\
Observe that theorem 6.3 gives the exact parameters (i.e. the minimum distance) of the codes $C_2(\mathcal{X})$ which cannot be found by F. Rodier in  \lbrack14, pp.207-208\rbrack.
\begin{theorem}
The minimum weight codewords correspond to quadrics which are pair of hyperplanes non-tangent to $\mathcal{X}$ such that the plane of intersection of the two hyperplanes intersects $\mathcal{X}$ at a non-singular hermitian plane curve.
\end{theorem}
\begin{theorem}
The second weight of the code $C_2(\mathcal{X})$ is $t^7-t^5-t^3$.\\
The codewords of second weight correspond to quadrics which are pair of hyperplanes non-tangent to $\mathcal{X}$ and the plane of intersection of the two hyperplanes intersecting $\mathcal{X}$ at a singular hermitian plane curve of rank 2 (i.e. a set of $t+1$ lines through a common point).
\end{theorem}
\begin{theorem}
The third weight of the code $C_2(\mathcal{X})$ is $t^7-t^5-t^2$.\\
The codewords of third weight correspond to quadrics which are pair of hyperplanes, one tangent to $\mathcal{X}$, the second non-tangent to $\mathcal{X}$ and such that the plane of intersection of the two hyperplanes intersects $\mathcal{X}$  at a non-singular hermitian plane curve.
\end{theorem}
\begin{theorem}
The fourth weight of the code $C_2(\mathcal{X})$ is $t^7-t^5$.\\
The codewords of fourth weight correspond to:\\
\textendash quadrics which are pair of hyperplanes, one tangent to $\mathcal{X}$, the second non-tangent to $\mathcal{X}$ and the plane of intersection of the two hyperplanes is tangent at the non-singular hermitian surface obtained from the second hyperplane.\\
\textendash quadrics which are pair of tangent hyperplanes 
and the plane of intersection of the two hyperplanes intersecting $\mathcal{X}$ at a single line.
\end{theorem}
\begin{theorem}
The fifth weight (for $t\ne 2,3$) of the code $C_2(\mathcal{X})$ is $t^7-t^5+t^3-t^2$.\\
The codewords of fifth weight (for $t\ne 2,3$) correspond to quadrics which are pair of tangent hyperplanes 
and the plane of intersection of the two hyperplanes intersecting $\mathcal{X}$ at a non-singular hermitian plane curve.
\end{theorem}
\section{Conjectures}
\subsection{Conjecture on the minimum distance of the code $C_h(\mathcal{X})$ in  $\mathbb{P}^{4}(\mathbb{F}_q)$} 
The author has tried to generalize the study of the code $C_2(\mathcal{X})$ to $C_h(\mathcal{X})$    where $\mathcal{X}$ is the non-degenerate hermitian variety defined by $x_{0}^{t+1}+x_{1}^{t+1}+ x_{2}^{t+1}+x_{3}^{t+1} +x_{4}^{t+1}= 0$ in $\mathbb{P}^{4}(\mathbb{F}_q)$ ($q=t^2$) and conjecture that: 
\paragraph{Conjecture 1}
\textit{
$$\mathrm{For}\quad h\le t\quad\# \mathcal{X}_{ Z(f)}(\mathbb{F}_{q} ) \le h( t^{5}+ t^{2})+t^3+1.$$
This upper bound is the best possible.
The minimum weight codewords correspond to hypersufaces which have the following configuration:\\
 \textendash hypersurfaces reaching the Serre-S\o rensen's upper bound for hypersurfaces (i.e. union of $h$ hyperplanes $H_i$ passing through a common plane $\mathcal{P}$)\\
 \textendash each hyperplane $H_i$ is non-tangent to $\mathcal{X}$\\
 \textendash and the plane $\mathcal{P}$ of intersection of the $h$ hyperplanes intersecting $\mathcal{X}$ at a non-singular hermitian plane curve.}\\\\
Unfortunately no proof has been found yet. In general, it seems that when we use this strategy of code-construction, the combinatorical complexity of finding the minimum distance increases drastically with the degree $h$.
\begin{Remark}
The preceding conjecture is true for $h=1$ (theorem \ref{BoseChar P4}) and for $h=2$ (theorems 6.1, 6.4) .
 \end{Remark} 
\subsection{Conjecture on the first five weights of the code $C_2(\mathcal{X})$ in  $\mathbb{P}^{n}(\mathbb{F}_q)$} 
The author has also tried to generalize the study of the code $C_2(\mathcal{X})$ to $\mathcal{X}$ the non-degenerate hermitian variety defined by $x_{0}^{t+1}+x_{1}^{t+1}+...+x_{N}^{t+1}= 0$ in $\mathbb{P}^{N}(\mathbb{F}_q)$ ($q=t^2$) and conjecture that:\\\\
\textbf{Conjecture 2}
\textit{
\begin{itemize}
\item[{1}] If $w_{i}$ $(1\le i\le 5)$ are the first five weights of the code $C_2(\mathcal{X})$, then there exist degenerate quadrics $\mathcal{Q}$ reaching the Serre-S\o rensen's upper bound for hypersurfaces (i.e. $\mathcal{Q}$ is a pair of distinct hyperplanes  $\mathcal{Q}=H_1\cup H_2$), giving codewords of weight $w_i$.
\item[{2}]The mimimum weight (i.e. $w_1$) codewords are only given by degenerate quadrics which are pair of distinct hyperplanes ($\mathcal{Q}=H_1\cup H_2$) such that $(H_1\cap H_2)\cap \mathcal{X}$ is a non-singular hermitian variety in $\mathbb{P}^{N-2}(\mathbb{F}_q)$ and: \\
 \textendash If N is even, the two hyperplanes $H_1$ and $H_2$ are non-tangent to $\mathcal{X}$\\
 \textendash If N is odd, the two hyperplanes $H_1$ and $H_2$ are tangent to $\mathcal{X}$.
\item[{3}] For $i>5$, there is no quadric which is a pair of distinct hyperplanes, giving codewords of weight $w_{i}$.
\end{itemize}}
Unfortunately no proof has been found yet.
\begin{Remark}
The second part of conjecture 2 is true for $N=3$ (see \lbrack 6, $\S$ 4.1-4.2\rbrack). The conjecture is also true for $N=4$: theorem 6.1, theorems 6.3 to 6.8 give the result. 
 \end{Remark}
\textbf{Acknowlegment}
The author thanks Prof. F. Rodier. His remarks and patience encouraged him to work on the problem.\\ \\
\textbf{References}\\
{\footnotesize \lbrack 1\rbrack   \ Y. Aubry,  Reed-Muller codes associated to projective algebraic varieties. In ``Algebraic Geomertry and Coding Theory ". (Luminy, France, June 17-21, 1991). Lecture Notes in Math., Vol. 1518, Springer-Verlag, Berlin, (1992), 4-17.\\
\lbrack 2\rbrack \ M. Boguslavsky, On the number of solutions of polynomial systems. Finite Fields and Their Applications 3 (1997), 287-299.\\ 
\lbrack 3\rbrack  \ R. C. Bose and I. M. Chakravarti, Hermitian varieties in finite projective space $PG(N,q)$. Canadian J. of Math.18 (1966), 1161-1182.\\
\lbrack 4\rbrack \  I. M. Chakravarti, Some properties and applications of Hermitian varieties in finite projective space PG(N,$q^2$) in the construction of strongly regular graphs (two-class association schemes) and block designs, Journal of Comb. Theory, Series B, 11(3) (1971), 268-283.\\
\lbrack 5\rbrack  \ F. A. B. Edoukou,  Codes defined by forms of degree 2 on quadric surfaces. 
arXiv:math. AG/0511679 v 1 28 Nov 2005. 5 pp., Submitted to IEEE. Trans.  Inform.  Theory  (2005).\\ 
\lbrack 6\rbrack  \ F. A. B. Edoukou, The weight distribution of the functional codes defined by forms of degree 2 on hermitian surfaces.  arXiv:math. AG/0512476 v 1 20 Dec 2005. 13 pp., Submitted to Finite Fields and Their Applications (2005).\\                     
\lbrack 7\rbrack  \ J. W. P. Hirschfeld, Projective Geometries Over Finite Fields (Second Edition) Clarendon  Press. Oxford 1998.\\
\lbrack 8\rbrack  \ J. W. P. Hirschfeld, Finite projective spaces of three dimensions, Clarendon press. Oxford 1985. \\
\lbrack 9\rbrack  \ J. W. P. Hirschfeld,  General Galois Geometies, Clarendon press. Oxford 1991.\\
\lbrack 10\rbrack  \ G. Lachaud, Number of points of plane sections and linear codes defined on algebraic varieties;  in " Arithmetic, Geometry, and Coding Theory ". (Luminy, France, 1993), Walter de Gruyter, Berlin-New York, (1996), 77-104. \\
\lbrack 11\rbrack \ D. B. Leep and L. M. Schueller, Zeros of a pair of quadric forms defined over finite field. Finite Fields and Their Applications 5 (1999), 157-176.\\ 
\lbrack12\rbrack \ E. J. F. Primrose,  Quadrics in finite geometries, Proc. Camb. Phil. Soc., 47 (1951), 299-304.\\
\lbrack13\rbrack \ D. K. Ray-Chaudhuri,  Some results on quadrics in finite projective geometrie based on Galois fields, Can. J. Math. Vol. 14 (1962), 129-138.\\
\lbrack14\rbrack \  F. Rodier, Codes from flag varieties over a finite field, Journal of Pure and Applied Algebra 178 (2003), 203-214.\\
\lbrack15\rbrack \  P. Samuel, G\'eom\'etrie projective, Presses Universitaires de France, 1986.\\
\lbrack 16\rbrack  \ J. P. Serre, Lettre \`a M. Tsfasman, In "Journ\'ees Arithm\'etiques de Luminy (1989)", Ast\'erisque 198-199-200 (1991), 3511-353.\\
\lbrack 17\rbrack  \ A. B. S\o rensen, Rational points on hypersurfaces, Reed-Muller codes and algebraic-geometric codes. Ph. D. Thesis, Aarhus, Denmark, 1991.\\
\lbrack 18\rbrack  \ H. P. F. Swinnerton-Dyer, Rational zeros of two quadratics forms, Acta Arith. 9  (1964), 261-270.\\
\lbrack 19\rbrack  \ J. Wolfmann, Codes projectifs \`a deux ou trois poids associ\'es aux hyperquadriques d'une g\'eom\'etrie finie, Discrete Mathematics 13 (1975), 185-211.}
\end{document}